\documentclass[twoside]{cleanMatt}
\author{Matthew Kvalheim\footnote{EECS Dept., University of Michigan, Ann Arbor, MI, USA (\texttt{kvalheim{@}umich.edu})} \qquad Shai Revzen\footnote{Depts. of EECS \& EEB, University of Michigan, Ann Arbor, MI, USA \texttt{shrevzen{@}umich.edu}}}
\journalname{}
\title{Reverse-engineering invariant manifolds with asymptotic phase}
\usepackage[utf8x]{inputenc}
\usepackage{thmtools}
\usepackage[T1]{fontenc}
\usepackage{lmodern}
\usepackage{amsmath}
\usepackage{amssymb}
\usepackage{amsthm}
\usepackage{tikz-cd}
\usepackage{subfig}
\newcommand{\pageFooter}{{\large Kvalheim \& Revzen}}
\newcommand{\concept}[1]{``\textit{#1}''}
\newcommand{\R}{\mathbb{R}}
\newcommand{\dist}[2]{\mathrm{dist}\left[{#1},{#2}\right]}
\newcommand{\bo}{\mathcal{O}}
\newcommand{\T}{\mathsf{T}}
\newcommand{\Nor}{\mathsf{N}}
\newcommand{\Hor}{\mathsf{H}}
\newcommand{\Ver}{\mathsf{V}}
\newcommand{\D}{\mathsf{D}}

\newcommand{\C}{\mathcal{C}}
\newcommand{\Q}{\mathcal{Q}}
\newcommand{\M}{\mathcal{M}}
\newcommand{\N}{\mathbb{N}}

\newcommand{\B}{\mathcal{B}}

\newcommand{\ct}{\setminus}

\newtheorem{Lem}{Lemma}

\newtheorem{Th}{Theorem}
\newtheorem{Co}{Corollary}

\newtheorem{Prop}{Proposition}

\newtheorem{prop-hand}{Proposition}[]
\newenvironment{Prop-hand}[1]
{\setcounter{prop-hand}{\numexpr#1-1\relax}\begin{prop-hand}}
{\end{prop-hand}}

\newtheorem{thm-hand}{Theorem}[]
\newenvironment{Th-hand}[1]
{\setcounter{thm-hand}{\numexpr#1-1\relax}\begin{thm-hand}}
	{\end{thm-hand}}

\theoremstyle{definition}
\newtheorem{Def}{Definition}
\newtheorem*{Def*}{Definition}

\theoremstyle{remark}
\newtheorem{Rem}{Remark}

\usepackage[numbers,sort&compress]{natbib}
\bibpunct{[}{]}{;}{a}{,~}{,}

\newcommand{\NOcitep}[2][]{}
\newcommand{\NOcitet}[2][]{}

\usepackage{marvosym}
\usepackage[%hyperref
    hypertexnames,%
    citecolor=blue,%
    colorlinks=true,%
    linkcolor=red%
]{hyperref}
\usepackage[all]{hypcap}

\begin{document}
\sffamily
\maketitle
\tableofcontents

% Put my document here

\section{Introduction}\label{sec:intro}
Most physical phenomena evolve through elaborate and difficult to model dynamical equations. 
The spectacular success of the physical sciences can be attributed to the fact that in some cases, simple models exist -- complex interactions lead to simpler solutions.
Our particular interest comes from biomechanics, where the notion of \concept{templates} and \concept{anchors} \citep{full_templates_1999} has been in use.
A template and anchor are a pair of models such that the template describes the essential features of a biomechanical behavior -- e.g., when running, animals bounce as if the center of mass is on a pogo stick \citep{blickhan1989spring} -- whereas the anchor is a more complete model that often contains specifics of the individual animal morphology.
One particular way to formalize a mathematical relationship between a template and anchor is to require \concept{asymptotic equivalence}, whereby the infinite horizon prediction of the template coalesces with that of the anchor.
This convergence onto a specific trajectory of the template is referred to in the dynamics literature as the \concept{asymptotic phase property} \citep{fenichel1973asymptotic,fenichel1977asymptotic,odeHale,smoothInvariant}, and is usually associated with \concept{normal hyperbolicity} of the template.
Numerous templates have been proposed by biologists, and these have fueled a flurry of activity in robotics whereby engineers have tried to design robots to express the self-same templates in their dynamics.
When successful, such attempts have enabled robots to run efficiently with legs \citep{galloway2010passive}, dynamically climb walls with little to no sensing \citep{lynch-2011-phd}, and reorient in free-fall to land safely \citep{Libby-2012-nature}.
In each of these cases, a specifically customized approach was used to anchor the template into the robot. 

Our contribution below is a general method for anchoring normally hyperbolic templates; furthermore, this method is universal in the sense that a broad class of normally hyperbolic template-anchor systems can be seen as special cases of our construction.
To express this formally, let $(\B,f)$ denote a smooth dynamical system consisting of a smooth manifold $\B$ and a vector field $f$ on $\B$.
Let $\phi^f_t(\cdot)$ denote the flow of $f$.
Suppose that $(\tilde \M, \tilde g)$ is another smooth dynamical system with flow $\phi^{\tilde g}_t(\cdot)$.
For simplicity, assume that both vector fields are \concept{complete} in the sense that their respective flows are respectively defined on all of $\B \times \R$ and $\tilde \M \times \R$.
One reasonable criterion for approximation is that there exists a  surjective map $\tilde P:\B\to\tilde\M$ which plays the role of \concept{model reduction} and satisfies $\tilde P\circ \phi_t^f = \phi_t^{\tilde g} \circ \tilde P$.
In other words, flowing a point $x$ forward by $\phi_t^f$ followed by reduction by $\tilde P$ is the same as flowing $\tilde P(x)$ forward by $\phi_t^{\tilde g}$.
We say that $\tilde P$ is a \concept{semi-conjugacy} between the families $\phi^f_t$ and $\phi^{\tilde{g}}_t$. 
If $\tilde P$ is at least $\C^1$, then the following two diagrams commute and are equivalent\footnote{To see that commutativity of the left diagram implies commutativity of the right, simply take partial derivatives of the left diagram with respect to $t$.
To show the opposite implication, just note that $\D \tilde P(f(x)) = \tilde g(\tilde P(x))$ implies that both $\tilde P\circ \phi^f_t(x)$ and $\phi^{\tilde g}_t\circ \tilde P(x)$ are the unique solution of the initial value problem $\dot y = \tilde g(\tilde P(y)),\,y(0)=\tilde P(x)$.}. 

\begin{equation}\label{diag:flow_invariant}
\begin{tikzcd}
\B \arrow{r}{\phi^f_t} \arrow{d}{\tilde P}
&\B \arrow{d}{\tilde P}
\\
\tilde \M \arrow{r}{\phi^{\tilde g}_t} &\tilde\M.
\end{tikzcd}
\qquad
\begin{tikzcd}
\B \arrow{r}{f} \arrow{d}{\tilde P}
&\T\B \arrow{d}{\D \tilde P}
\\
\tilde \M \arrow{r}{\tilde g} &\T \tilde M.
\end{tikzcd}
\end{equation}
If $\tilde P$ is a smooth submersion, we also say that $f$ is a \concept{lift} of $\tilde g$.

Given a dynamical system $(\B,f)$ and a surjective submersion $\tilde P:\B\to \tilde \M$ with connected fibers (level sets),
it is a textbook exercise to show that $f$ is a lift of some vector field $\tilde g$ on $\tilde \M$ if and only if the Lie bracket $[h,f] \in \ker \D \tilde P$ whenever $h$ is a vector field taking values in $\ker \D \tilde P$ (\citet{lee2013smooth} Exercise 8-18).
% * <shrevzen@umich.edu> 2016-08-15T23:03:41.004Z:
%
% The notation [h,f] hasn't been defined. Is this a Lie bracket?
%
% ^.

It is similarly straightforward to show that given any dynamical system $(\tilde \M,\tilde g)$ and smooth submersion $\tilde P:\B\to\tilde \M$, there exist many vector fields $f$ on $\B$ which are lifts of $\tilde g$ (\citet{lee2013smooth} Exercise 8-18).
It is easy to give explicit formulas for such an $f$; for example, one may take $f(x):= \D \tilde P_x^\dagger \tilde g(\tilde P(x))$, where $(\cdot)^\dagger$ is the Moore-Penrose pseudoinverse of a linear map discussed in Appendix \ref{app:pseudoinverse}.

A less trivial problem concerns the case in which $(\tilde \M,\tilde g)$ is embedded in $\B$ in the sense that there is an embedding $F:\tilde \M\to \B$ so that $\M:=F(\tilde \M)$ is a properly embedded submanifold of $\B$ and so that the vector field $g:\M \to \T \B$ defined by $g:=\D F \tilde g$ is everywhere tangent to $\M$.
Given a surjective smooth submersion $\tilde P:\B\to \tilde\M$, we define the map $P:\B\to\B$ by $P:= F \circ \tilde P$.
Assuming that $\tilde P|_{\M} = F^{-1}$, it follows that $P$, when thought of as a map into $\M$, is a submersion and retraction (i.e., $P(P(x))= P(x)$ for every $x \in \B$).
The goal is to produce a vector field $f:\B\to\T \B$ such that $f|_\M = g$, $f$ is a lift of $g$ (in the sense that $\forall x \in \B: \D Pf(x) = g(P(x))$), and $\M$ is an asymptotically stable invariant manifold of $\B$.
In this case, the na\"ive solution of the preceding paragraph will no longer work, as 
$\D P_x ^\dagger g(P(x))$ will generally no longer by tangent to $\M$ for $x \in \M$, so $\M$ will not be an invariant manifold of $f$.
The main contribution of our work (in \S \ref{sec:construction}) is to give an explicit formula for constructing a vector field $f$ with flow $\phi_t(\cdot)$ such that $f$ is a lift of $g$ given $\B,\tilde\M,\tilde g,F$, and a couple other ingredients, at least in the case that the embedding $\M$ is the level set of a smooth submersion; i.e., if $\M$ can be expressed as the set of points satisfying some constraints (which must be full rank).
In the case we consider, in which $\B$ is an open subset of $\R^n$, this formula for $f$ will involve only standard matrix computations.
Furthermore, $f$ will render $\M$ asymptotically stable and normally hyperbolic for suitable parameter choices.

Our construction works by separately constructing the components of a vector field $f:= f_h + f_v$, where we refer to $f_h$ as the \concept{horizontal} vector field and to $f_v$ as the \concept{vertical} vector field $f_v$.
The vertical vector field is responsible for the stabilization of $\M$, and the horizontal vector field is responsible for ensuring that $f$ is a lift of $g$.

Since $\M$ is a subset of $\B$, different language can be used to discuss model reduction and there is further justification for the value of $P$ as a model reduction tool.
Following \citet{fenichel1973asymptotic,fenichel1977asymptotic,odeHale,smoothInvariant}, we make the following definition:

	\begin{Def}\label{def:asymptotic_phase}
		Suppose that $\M$ is asymptotically stable under the flow $\phi_t(\cdot)$ with basin of attraction $\B$.
		We say that $\M$ has \concept{asymptotic phase} if there exists a map $P:\B \to \B$ with $P(\B)=\M$ such that for any $x \in \B$: $$\lim_{t\to\infty}\|\phi_t(x)-\phi_t(P(x))\| = 0.$$
		We refer to $P$ as the \concept{phase map} or simply as \concept{phase}, and say that $\M$ has $\C^k$ asymptotic phase if the map $P:\B\to \B$ is $\C^k$.
	\end{Def}

If $\M$ has asymptotic phase as the result of our construction, then not only do trajectories in $\B$ approach $\M$; they approach \textit{specific trajectories} in $\M$.
This is the most important reason that our approach yields a dynamical system on $\B$ for which the dynamics on $\M$ are a good approximation.
Unlike approximation techniques like linearization, approximation using our map $P$ on $\B$ does not get worse on longer time intervals; if the dynamics are deterministic, \textit{the approximation gets better} the longer the time of execution is.

The remainder of this paper is organized as follows.
In \S \ref{sec:motivation}, we first motivate our construction by proving Theorem \ref{th:main_theoretical}, which shows that the vector field $f$ on the basin of attraction of a normally hyperbolic invariant manifold always admits a certain decomposition $f = f_h + f_v$.
In \S \ref{sec:construction}, we then describe our construction, which works by explicitly constructing vector fields $f_h$ and $f_v$, in the case that $\M$ is the regular level set of a submersion $G:\B \to \R^{n-k}$. 
In \S \ref{sec:our_persistence}, we prove that our construction renders $\M$ normally hyperbolic with asymptotic phase, and both $\M$ and its asymptotic phase \concept{persist} under perturbations of our constructed vector field.
In \S \ref{sec:examples} we give a few simple, but hopefully illustrative, examples on how our results can be applied to physical systems.
In \S \ref{sec:lyap_extend}, we show that Lyapunov functions for the dynamics on $\M$ extend in a natural way to Lyapunov functions on all of $\B$, the basin of attraction of $\M$ under our constructed dynamics.
In \S \ref{sec:generalizing} we discuss topological constraints arising from our assumptions and indicate that our approach might generalize to the case in which $\M$ is not a regular level set.
The paper concludes with suggestions for future work and discussion of the relevance of this work to physical and biological systems of practical interest.

We assume basic familiarity with linear algebra, smooth dynamical systems theory, and point-set and differential topology.
In the appendices, we briefly review relevant notions from these different subjects.
We also give a brief overview of definitions and results from the theory of normally hyperbolic invariant manifolds.
In this work, we make some use of notions from the mathematical framework of \concept{fibered manifolds}, \concept{fiber bundles}, and \concept{connections}; see, e.g., \citet{kolar1999natural}.
In an attempt to make our work more self-contained, we briefly summarize relevant aspects of these concepts in Appendix \ref{app:fiber_connections}, and prove some results on connections in Appendix \ref{app:connection_lemmas} which we will have use for in proving some of our main results.
%Theorem \ref{th:main_theoretical} in \S \ref{sec:motivation}.

\section{Motivation: decomposition of NHIM-defining vector fields}\label{sec:motivation}
Normally hyperbolic invariant manifolds (NHIMs) are generalizations of hyperbolic fixed points and periodic orbits.
NHIMs have desirable properties from the perspective of model reduction, as will be discussed in \S \ref{sec:NHIM_implications}.
We will now discuss a decomposition of the vector field $f$ on the stability basin $\B$ of a $k$-dimensional NHIM $\M$ which lends insight into how one might \textit{construct} dynamics on an open set $\B$ containing $\M$ rendering $\M$ exponentially stable and normally hyperbolic in the case that asymptotic phase is at least $\C^2$.

The key idea behind the following Theorem \ref{th:main_theoretical} is that $f$ may be decomposed into two parts playing distinct roles: $f_v$ ensures that $\M$ is asymptotically stable, and $f_h$ ensures that $P$ is the asymptotic phase of $\M$.
For stable NHIMs additionally possessing $\C^2$ asymptotic phase, such a natural decomposition of the flow in the stability basin is always possible, and inverting it provides a design tool for specifying anchor dynamics that produce a desired template dynamic. 
The idea behind this result produced the construction in \S \ref{sec:construction} which, under certain assumptions, allows one to take a vector field on a compact submanifold $\M$ of an open set $\B \subseteq \R^n$ together with a phase map $P:\B \to \B$ and explicitly construct a vector field on $\B$ such that $\M$ is a NHIM.
Fiber bundles and connections, used below, are defined in Appendix \ref{app:fiber_connections}.

\begin{Th}\label{th:main_theoretical}
	Let $r \geq 2$ and let $\M$ be a $\C^r$ exponentially stable NHIM of the $\C^r$ vector field $f$ with basin of attraction $\B \subseteq \R^n$ and additionally possessing
	$\C^r$ asymptotic phase $P:\B \to \B$.
	Then there exists a decomposition of $f = f_h +f_v$, with $f_v$ taking values in $\Ver \B := \ker \D P$ and with $\D P f_h(x) = f(P(x))$ for all $x \in \B$. 
	In particular, $f_h|_\M = f|_\M$.
\end{Th}

\begin{proof}
	As we will see in Proposition \ref{prop:phase}, $(B,P,\M)$ is a $\C^r$ fibered manifold.
	Lemma \ref{lem:connection} in Appendix \ref{app:connection_lemmas} shows that there exists a $\C^{r-1}$ connection $\Hor \B$ for $P$ on $\B$ which restricts to $\T \M$ on $\M$.
	Letting $\Ver \B:= \ker \D P$, we see that $\B = \Ver \B \oplus \Hor \B$, 
	and we may uniquely write $f = f_h + f_v$, with $f_h$ taking values in $\Hor \B$ and $f_v$ taking values in $\Ver \B$.
	Lemma \ref{lem:hor_lift_smooth} in Appendix \ref{app:connection_lemmas} shows that $f_h$ is $\C^{r-1}$, and hence $f_v := f - f_h$ is also $\C^{r-1}$.
\end{proof}

Next, we will informally state some results from the theory of NHIMs which we will use to prove Theorem \ref{th:persistence_of_our} in \S \ref{sec:our_persistence}.
First, in \S \ref{sec:NHIM_structure}, we state Proposition \ref{prop:phase} describing the structure of the basin of attraction of a NHIM -- in particular, this proposition shows that asymptotically stable NHIMs have unique asymptotic phase.
Next, in \S \ref{sec:NHIM_persistence}, we state Propositions \ref{prop:persistence} and \ref{prop:converse_mane} -- Proposition \ref{prop:persistence} shows that NHIMs and the structure of their stability basins persist under small perturbations of their defining vector fields, and Proposition \ref{prop:converse_mane} shows that \textit{every} invariant manifold persisting under small perturbations of its defining vector field is normally hyperbolic.
Finally, in \S \ref{sec:NHIM_implications} we discuss the implications of these results for model reduction. 

We will consider only compact, asymptotically stable NHIMs in what follows.
All NHIMs we consider are embedded submanifolds of $\R^n$, and the same is true for all NHIMs in Propositions \ref{prop:phase}, \ref{prop:persistence}, and \ref{prop:converse_mane}.
Precise definitions and statements of the following results are given in Appendix \ref{app:NHIMs}.

\subsection{The structure of the basin of attraction of a NHIM}\label{sec:NHIM_structure}
We defined \concept{asymptotic phase} in Definition \ref{def:asymptotic_phase}.
$x\in \B$ has the same asymptotic phase as $q \in \M$ if $\phi_t(x)$ and $\phi_t(q)$ asymptotically coalesce.
However, it could conceivably be the case that for some point $x \in \B$, there are multiple points $q_1,q_2\in\M$ such that $\phi_t(x)$ asymptotically coalesces with both $\phi_t(q_1)$ and $\phi_t(q_2)$.
A more stringent concept that removes this ambiguity is that of \concept{unique asymptotic phase}. 
If $\M$ has unique asymptotic phase, it can still be the case that for some $x \in \B$, $\phi_t(x)$ can asymptotically coalesce with multiple points of $\M$ -- however, there will be a \textit{unique} point of $\M$ that $\phi_t(x)$ asymptotically coalesces with \textit{most rapidly}.
We will use this terminology in stating results in this section.
Following (\citet{odeHale} p. 217) and \S I.E. of \citet{fenichel1973asymptotic}, we give the following definition making the notion of \concept{unique asymptotic phase} precise.

\begin{Def}\label{def:unique_asymptotic_phase}
	We say that $\M$ has \concept{unique asymptotic phase} if $\M$ has asymptotic phase $P:\B\to\B$ and additionally for any $x \in \B$ and any $q \in \M$ not equal to $P(x)$, 
	$$\lim_{t\to\infty}\frac{\|\phi_t(x)-\phi_t(P(x))\|}{\|\phi_t(x)-\phi_t(q)\|}=0.$$
	We say that $\M$ has \concept{$\C^k$ unique asymptotic phase} if $\M$ has unique asymptotic phase and if the map $P:\B\to \B$ is $\C^k$.
\end{Def}

Next, we use this definition in stating Proposition \ref{prop:phase} below which is a combination of results from \citet{fenichel1973asymptotic,fenichel1977asymptotic}, and Theorem 4.1 of \citet{hirsch1977}.
Informally, an asymptotically stable invariant manifold $\M$ is \concept{$r$-normally hyperbolic} if, to first order, trajectories of $f$ approach $\M$ $r$-times faster than nearby trajectories in $\M$ converge together in positive time.  
Without further qualification, \concept{normally hyperbolic} and \concept{NHIM} refers to a NHIM which is at least $1$-normally hyperbolic.
A precise statement of this proposition is given in Appendix \ref{app:NHIMs}.

\begin{Prop}\label{prop:phase}
	Let $\M \subseteq \R^n$ be a compact $\C^r$ $k$-dimensional asymptotically stable NHIM, invariant under the flow $\phi_t(\cdot)$ of the vector field $f:\Q\to \T \Q$ defined on an open neighborhood $\Q \subset \R^n$ of $\M$.
	Then the following holds:

	\begin{enumerate}
		\item The stability basin $\B$ of $\M$ is partitioned into codimension-$k$ $\C^r$ manifolds $(W_q)_{q\in\M}$ permuted under the flow.
		Explicitly, $\phi_t(W_q) = W_{\phi_t(q)}$.
		Each $W_q$ is $\C^r$ diffeomorphic to $\R^{n-k}$.
		Each $W_q$ intersects $\M$ transversally in the point $q$.
		
		\item Let $P:\B \to \B$ be the map that sends $x \in \B$ to $q$, where $x \in W_q$. 
		Then $P$ is a continuous map, and $(\B,P,\M)$ is a $\C^0$ fibered manifold with $(n-k)$-dimensional Euclidean fibers.
		
		\item $\M$ has unique asymptotic phase given by $P:\B\to \B$. 				
		
		\item 
		If the flow satisfies an additional condition related to the rate of expansion of the flow on $\M$\footnote{See the precise statement of this proposition in Appendix \ref{app:NHIMs}.}, then the phase map $P$ is $\C^{r-1}$. It additionally follows that $(\B,P,\M)$ is a $\C^{r-1}$ fibered manifold with $(n-k)$-dimensional Euclidean fibers.
	
	\end{enumerate}	
	
\end{Prop}

\subsection{Persistence of this structure under perturbations}\label{sec:NHIM_persistence}

Informally, two differentiable functions are \concept{$\C^1$-close} if both the functions and their derivatives are close at all points of their domains.
Two embeddings of a given manifold are \concept{$\C^1$-close} if the functions defining the embeddings are $\C^1$-close.
A precise definition is given in Definition \ref{def:C1_close} in Appendix \ref{app:NHIMs}.
We will use this concept in stating Proposition \ref{prop:persistence}, a more precise statement of which is also given in Appendix \ref{app:NHIMs}.
Proposition \ref{prop:persistence} is a robustness result; it gives conditions under which $\M$ and its unique asymptotic phase persist under $\C^1$-small perturbations by $\C^r$ vector fields.

\begin{Prop}\label{prop:persistence}
	Let $\M \subseteq \R^n$ be a compact $\C^r$ $k$-dimensional asymptotically stable NHIM, invariant under the flow $\phi_t(\cdot)$ of the vector field $f:\Q\to \T \Q$ defined on an open neighborhood $\Q \subset \R^n$ of $\M$.
	Then the following holds:
	
	\begin{enumerate}
		\item Let $g:\Q \to \T \Q$ be another $\C^r$ vector field which is sufficiently $\C^1$-close to $f$.
		Then there is a unique $\C^r$ embedded submanifold $\M'$, $\C^r$ diffeomorphic to $\M$, $\C^1$-close to $\M$, and a NHIM for the vector field $g$.
		Furthermore, the fibers $W_q$ persist; i.e., there is a unique partition of the  stability basin $\B'$ of $\M'$ into codimension-$k$ $\C^r$ manifolds $W_{q'}'$ satisfying all of the properties with respect to $g$ and $\M'$ which were satisfied by the manifolds $W_p$ with respect to $f$ and $\M$.
		The fibers $W_{q'}'$ are $\C^1$-close to those of $W_q$ on $\B \cap \B'$.
		$\M'$ has unique asymptotic phase $P'$ whose fibers are $W_{q'}'$, and $P':\B' \to \B'$ is a continuous function.
		
		\item 
		If the flow satisfies an additional condition related to the rate of expansion of the flow on $\M$\footnote{See the precise statement of this proposition in Appendix \ref{app:NHIMs}.}, then the phase map $P'$ corresponding to the perturbed vector field $g$ is also $\C^{r-1}$ if $g$ is sufficiently $\C^1$-close to $f$, and $(\B',P',\M')$ is a $\C^{r-1}$ fibered manifold with $(n-k)$-dimensional Euclidean fibers.
		Under these conditions, $\M'$ has unique $\C^{r-1}$ asymptotic phase $P'$. 		
	\end{enumerate}
\end{Prop}

Next, without being too precise, we state a result of \citet{mane1978persistent} stating a partial converse result.

\begin{Prop}\label{prop:converse_mane}
	Let $\M$ be a compact invariant manifold of the $\C^1$ vector field $f$ which persists under $\C^1$-small perturbations to $f$.
	Then $\M$ is normally hyperbolic.
\end{Prop}

\subsection{Implications of NHIM results for model reduction}\label{sec:NHIM_implications}
Propositions \ref{prop:phase} and \ref{prop:persistence} show that normally hyperbolic invariant manifolds have great utility in model reduction.

If a dynamical system possesses an asymptotically stable NHIM $\M$, Proposition \ref{prop:phase} says that not only do trajectories in the basin of attraction $\B$ approach $\M$, but these trajectories actually approach \textit{specific trajectories in $\M$}.
This means that trajectories in $\B$ can be approximated by trajectories in $\M$, justifying an approximation of the dynamics on $\B$ by the dynamics restricted to $\M$.
The unique asymptotic phase property implies that corresponding to each trajectory in $\B$, there is a \textit{unique} best approximating trajectory in $\M$.

Proposition \ref{prop:persistence} shows that normally hyperbolic invariant manifolds are robust; they persist under small perturbations of the vector field.
This is important from a physical modeling perspective -- since measurements of physical quantities can only be obtained with finite precision, persistence of NHIMs has the implication that they remain present in models despite small parameter errors, and thus can actually meaningfully represent features of the physical world.

Proposition \ref{prop:converse_mane} shows in a precise sense that \textit{every} compact invariant manifold persisting under all small perturbations is normally hyperbolic.
This means that every model reduction that produces a robust model on a compact manifold arises from the class of models we discuss here.

\section{Attractors arising as regular level sets}\label{sec:construction}
In this section, we present our main contribution: we \concept{reverse-engineer} the results of \S \ref{sec:motivation}, producing a general algorithm for anchoring templates in anchors in such a way that the template is an asymptotically stable NHIM.
In order to make our construction more concrete and to increase the accessibility of this work, we only consider Euclidean ambient spaces in what follows.
The reader fluent in Riemannian geometry will recognize that much, if not all, of our construction can be generalized to the case in which $\B$ (used below) is an open subset of a Riemannian manifold.
\subsection{The setup}
Let $\tilde \M$ be a compact, connected $k$-dimensional $\C^r$ manifold ($r \geq 1$) and $\tilde g$ be a $\C^r$ vector field on $\tilde \M$.
The $\tilde{\M}$ manifold is an abstracted template, allowing the template dynamics to be described in whatever representation is most natural.
Let $\B$ be a connected open subset of $\R^n$, equipped with the standard Euclidean inner product $\langle \cdot , \cdot \rangle$ and norm $\|\cdot\|$.  
Let $F: \tilde \M \hookrightarrow \B$ be a proper $\C^r$ embedding, and denote $\M := F(\tilde \M)$ and $g:= \D F (\tilde g)$.
Here $\M$ is the concrete instance of the template that appears in the anchor's state space.

In this section, we describe an approach to defining a vector field on $\B$ rendering $\M$ asymptotically stable with unique asymptotic phase and normally hyperbolic, under the primary assumption that $\M$ is the level set of a submersion $G:\B \to \R^{n-k}$.
We accomplish this by constructing a specific connection on $\B$ and then constructing vector fields $f_h$ and $f_v$ as in \S \ref{sec:motivation}.
We now state our assumptions in detail. 

\subsection{Assumptions}
\label{sec:assumptions}

\subsubsection{The attractor as a regular level set}\label{sec:level_set}
Let $G:\B \to \R^{n-k}$ be a $\C^r$ submersion.
This implies that $G^{-1}(0)$ is a $k$-dimensional $\C^r$ embedded submanifold of $\B$.
We require that $\M = G^{-1}(0)$, i.e., we know how to define the attractor using a set of simultaneous constraints which are never redundant, i.e. all constraints are always active.
We make the additional assumption that $\M$ is compact.
We write the $i$-th component of $G$ as $G_i$, so that 
\begin{equation}
G(x) = \left(G_1(x),\ldots,G_{n-k}(x)\right).
\end{equation}

\subsubsection{The phase map}\label{sec:phase_assump}
We assume there exists a $\C^r$ map $\tilde P:\B \to \tilde \M$ satisfying $\tilde P|_{\M} = F^{-1}$.
Define the map $P:\B \to \B$ by $P:=F\circ \tilde P$. 
It follows that $P$ fixes the set of points in $\M$.
Viewed as a map into $\M$, $P$ is a submersion and a retraction.
Note that for $x \in \M$, $P^{-1}(x) = \tilde P^{-1}(\tilde P(x))$.

\subsubsection{A specific connection}\label{sec:connection}
We assume 
\begin{equation}\label{eq:DG_conn}
\forall x \in \B:\, \T_x \B = \ker \D P_x \oplus \ker \D G_x,
\end{equation}
so that $\ker \D G$ is a connection (as defined in Appendix \ref{app:fiber_connections}) for $P$ on $\B$.
Note that $\ker \D P_x = \ker \D \tilde P_x$ since $\D P_x = \D F_{\tilde P(x)} \D \tilde P_x$ and $\D F_{\tilde P(x)}$ is full rank, so we could have written equation \eqref{eq:DG_conn} using $\ker \D \tilde P_x$.
Intuitively, this means any two states that coalesce asymptotically must be distinguishable to first order using the constraint functions.

\subsubsection{Completeness}\label{sec:completeness}
This last assumption will be used to ensure that vector fields we define later are complete, i.e. the trajectories of these vector fields never leave $\B$ and never tend to $\infty$ in norm in finite time. 
We assume that $\forall x \in \B: \|G(x)\| < \sup_{y\in\B}\|G(y)\|$ and that $\|G(x)\| \to \sup_{y\in\B}\|G(y)\|$ as $x$ approaches any point of $\partial \B \cup \{\infty\}$.
		
For later purposes, we define the function $V:\B\to \R$ by 
\begin{equation}
V(x):=\|G(x)\|^2 = \sum_{i=1}^{n-k}G^2_i(x),
\end{equation}
where $G_i:\B\to\R$ is the $i$th component function of $G$.

\subsection{Consequences of the assumptions}
\label{sec:consequences}
We define the family of projections $\Pi^P:\T\B\to\ker \D P$ for each $x \in \B$ by
\begin{equation}\label{eq:Pi_P_def}
\Pi^P_x:= [I - \D P_x^\dagger  \D P_x].
\end{equation}
where $(\cdot)^\dagger$ denotes the Moore-Penrose pseudoinverse relative to the standard inner product on $\R^n$. 	
We will occasionally suppress the subscript $x$ in the sequel when the notation becomes cumbersome, unless we wish to emphasize the role of $x$. 
Note that $\Pi_x^P:\T_x\B\to \ker \D P_x$.

\begin{Lem}\label{lem:VQ_props}
	The function $V$ satisfies:
	\begin{enumerate}
		\item $V^{-1}(0) = \M$.
		\item $\forall x \in \B \ct \M: V(x)> 0$.
		\item $\forall x \in \B \ct \M: \Pi^P_x\nabla V(x) \not = 0$,
	\end{enumerate}
	
\end{Lem}
Geometrically, the third property says that there always exist vectors tangent to the fibers of $P$ pointing in some direction along which the value of $V$ changes.
Note that the third condition implies that $V$ is a submersion on $\B \ct \M$, because $V:\B \to \R$ is a submersion if and only if $\nabla V(x)\not = 0$ for all $x\in\B$.

\begin{proof}
	Since $V$ is the square of the Euclidean norm of $G$, $V^{-1}(0) = G^{-1}(0) = \M$.
	Also, since $G$ is nonzero off of $\M$, the norm of $G$ is strictly positive off of $\M$ and hence $V$ is strictly positive on $\B \ct \M$.
	It remains only to show the third property.
	
	The fibers of $P$ are $(n-k)$-dimensional, the codomain of $G$ is $(n-k)$-dimensional, and the fibers of $P$ intersect the level sets of $G$ transversally (by equation \eqref{eq:DG_conn}).
	These facts together with the inverse function theorem imply that 
	the restriction of $G$ to any fiber of $P$ is a local diffeomorphism.
	If $x \not \in \M$, then $G(x) \not = 0$.
	Since the function $y \mapsto \|y\|^2$ is a submersion on $\R^{n-k} \ct \{0\}$, it follows that the restriction $V|_{P^{-1}(P(x))}:P^{-1}(P(x)) \to \R$ is a submersion at $x$ since the composition of submersions is again a submersion.
    Hence 
    \begin{equation*}
    \D V_x \Pi^P_x
    \end{equation*}
    is nonzero.
    For any $x \in \B$ and any $w \in \T_x\R^n$ with $\D V_x \Pi^P_x w \not = 0$, we have $\D V_x \Pi^P_x w = \langle \nabla V(x),\Pi^P_x w\rangle = \langle \Pi^P_x \nabla V(x),w\rangle$ (since orthogonal projections are self-adjoint), 
    so it follows that $\Pi^P_x \nabla V(x)\not = 0$.
\end{proof}

    \begin{Lem}\label{lem:compact_sublevel}
    	For $d < \sup_{x \in \B}V(x)$, every sublevel set $V^{-1}((-\infty,d])$ is compact.
    \end{Lem}
    \begin{proof}
    	Note that by composing $V$ with any diffeomorphism of $\R$ onto the interval $(-1,1)$ we may assume that $\sup_{x \in \B}V(x) < \infty$. 
    	Then $V$ extends to a continuous function $\tilde V: \bar \B \to \R$ by defining $\tilde V$ to be equal to $\sup_{x \in \B}V(x)$ everywhere on $\partial \B$.
    	It follows that if $d < \sup_{x \in \B}V(x)$, $V^{-1}((-\infty,d]) = \tilde V^{-1}((-\infty,d])$ is a closed subset of the closed subset $\bar \B$,
    	so $V^{-1}((-\infty,d])$ is closed in $\R^n$. 
    	Furthermore, $V^{-1}((-\infty,d])$ is bounded since the fact that $V(x) \to \sup_{x \in \B}V(x)$ as $x \to \infty$ implies that there is some bounded set outside of which $V(x) > d$, so $V^{-1}((-\infty,d])$ must be contained in this bounded set.
    	The Heine-Borel Theorem now implies that $V^{-1}((-\infty,d])$ is compact.
    \end{proof}
    
    \begin{Co}\label{co:compact_G_level}
    	Level sets of $G$ are compact.
    \end{Co}
    \begin{proof}
    	Given any $z\in\R^n$, $G^{-1}(z)$ is closed in $\B$ and hence also closed in the compact subset $V^{-1}(\|z\|^2)\subseteq \B$.
    \end{proof}

    \begin{Lem}\label{lem:complete_connection}
    	If $r\geq 2$, then $\Hor \B:= \ker \D G$ is a \concept{complete} connection (defined in Appendix \ref{app:fiber_connections}).
    \end{Lem}
    \begin{proof}
    	Let $\gamma:[0,1]\to \M$ be any $\C^r$ path.
    	Then $\dot \gamma$ defines a $\C^{r-1}$ vector field on the closed set $\gamma([0,1])$, which may be extended to a $\C^{r-1}$ vector field $g: \M \to \T \M$ on all of $\M$ (\citet{lee2013smooth} Chapter 8).
    	Lemma \ref{lem:hor_lift_smooth} shows that the horizontal lift $\tilde g$ of $g$ is a $\C^{r-1}$ vector field.
    	The horizontal lift of $\gamma$ with any initial point $q \in P^{-1}(\gamma(0))$ is the solution to the initial value problem $\dot x = \tilde g(x); x(0) = q$ on the time interval $[0,1]$.
    	The theory of ordinary differential equations ensures that a solution exists on some time interval (\citep{odeHirschSmale} p. 162); furthermore, if it can be guaranteed that the solution is confined to some fixed compact set, then the solution exists for all time (\citep{odeHirschSmale} page 171).
    	Since $\Hor \B = \ker \D G$, the solution of the stated initial value problem is confined to a single level set of $G$; by Corollary \ref{co:compact_G_level}, this level set is compact.
    	It follows that the solution of the initial value problem exists for all time, and hence the lift of $\gamma$ is defined on all of $[0,1]$. 
    	Since $\gamma$ was arbitrary, this completes the proof.
    \end{proof}

\begin{Prop}\label{prop:euclidean_fibers}
Assume $r \geq 2$. 
The fibers of $P$ are $\C^r$ diffeomorphic to $\R^{n-k}$
\end{Prop}
\begin{proof}
Consider the vector field $h:\B\to \T \B$ defined by $h(x):=-\Pi^P_x \nabla V(x)$.
As will be shown in the proof of Proposition \ref{prop:fv}, $\M$ is a globally asymptotically stable invariant manifold of $h$, and each fiber of $P$ is invariant under the flow of $h$.
Since $P|_\M$ is the identity, each fiber of $P$ intersects $\M$ in a single point.
It follows that the flow of $h$ restricted to any fiber of $P$ yields a well-defined flow on this fiber with a globally asymptotically stable equilibrium point.
Since $h$ is a $\C^{r-1}$ vector field, Theorem 2.2 of \citet{wilson1967structure} shows that each fiber of $P$ is $\C^{r-1}$ diffeomorphic to $\R^{n-k}$.
Theorem 2.10 on page 52 of \citet{hirsch1976differential} shows that two $\C^s$ manifolds are $\C^s$ diffeomorphic if and only if they are $\C^1$ diffeomorphic, where $s \geq 1$.
Since we assumed $r \geq 2$, this completes the proof. 
\end{proof}

\begin{Co}
If $r \geq 2$, $(\B,\tilde P, \tilde \M, \R^{n-k})$ is a $\C^r$ fiber bundle. 
\end{Co}
\begin{proof}
Lemma \ref{lem:complete_connection} showed that if $r \geq 2$, then $\Hor \B:= \ker \D G$ is a complete connection. 
It is a standard result that this implies that $\tilde P:\B \to \tilde \M$ defines a $\C^r$
% Strictly speaking, the standard proof in e.g. Del Hoyo only gives you $\C^{r-1}$; you then need to use Whitney approximation on paths to replace his $\C^{r-1}$ local trivializations with $\C^r$ ones.
fiber bundle\footnote{Though, interestingly, often proved incorrectly \citep{del2015complete}.} \citep{del2015complete,kolar1999natural}.
Proposition \ref{prop:euclidean_fibers} showed that the fibers are $\C^r$ diffeomorphic to $\R^{n-k}$, completing the proof.
\end{proof}

\subsection{Lifting the dynamics on $\M$ to $\B$}
We will define the vector field $f_h:\B \to \T \B$ to be the vector field such that
$\forall x \in \B,$ $f_h(x)$ is the unique vector in $\Hor_x \B$ such that $\D P_x f_h(x) = g(P(x))$.
$f_h$ is the \concept{horizontal lift} of $g$ as defined in Appendix \ref{app:fiber_connections}.
It is shown in Lemma \ref{lem:hor_lift_smooth} in Appendix \ref{app:connection_lemmas}
that $f_h$ exists, is unique, and is $\C^{r-1}$, but this will also be clear from the following construction.
Given $g: \M \to \T \M$, we now explicitly construct the lift $f_h$ of $g$ in a form amenable to matrix computation. 
We use the standard identifications of $\T_x\R^n$ with $\R^n$ and $\Ver_x\B, \Hor_x\B$ as subspaces of $\R^n$.
Given $x \in \B$, define the matrix-valued function $T: \B \to L(\R^n, \R^n)$ by 
\begin{equation}\label{eq:T_def}
T_x := \left[I - \D G_x^\dagger \D G_x\right]\D P_x^\dagger \D P_x + \D G_x^\dagger \D G_x.
\end{equation}

\begin{Lem}\label{lem:T_lem}
	For all $x\in \B$, $T_x$ is invertible. 
	$T_x$ maps $\ker \D G_x$ isomorphically onto itself and $T_x$ maps $\ker \D P_x$ isomorphically onto $(\ker \D G_x)^\perp$. 
	In particular, $T_x$ isomorphically maps $\ker \D G_x$ and $\ker \D P_x $ onto subspaces which are orthogonal with respect to the Euclidean inner product.
\end{Lem}

\begin{proof}
	Let $x \in \B$ and let $w \in \T_x \B$. 
	Equation \eqref{eq:DG_conn} tells us that $w = w_P + w_G$, with $w_P \in \ker \D P_x$ and $w_G \in \ker \D G _x$.
	Hence 
	\begin{equation}\label{eq:Tx_split}
	T_x w = \left[I - \D G_x^\dagger \D G_x\right] \D P_x^\dagger  \D P_x w_G + \D G_x^\dagger \D G_x w_P
	\end{equation}
	From equation \eqref{eq:Tx_split}, it follows that $T_x$ maps $\ker \D G_x$ into $\ker \D G_x$ and $T_x$ maps $\ker \D P_x$ maps into $(\ker \D G_x)^\perp$.
	Equation \eqref{eq:DG_conn} implies that $\ker \D G_x \cap \ker \D P_x = \emptyset$.
	Since $\ker \D G_x = \ker \D G_x^\dagger \D G_x$, this fact and the rank-nullity theorem imply that $\D G_x ^\dagger \D G_x$ maps $\ker \D P_x$ isomorphically onto $(\ker \D G_x)^\perp$ (since $\ker \D G_x$ and $\ker \D P_x$ have complimentary dimension).
	A similar argument shows that $\left[I - \D G_x^\dagger \D G_x\right] \D P_x^\dagger  \D P_x$ maps $\ker \D G _x$ isomorphically onto itself. 
	This completes the proof, since invertibility of $T_x$ follows from the algebraic fact that if a linear map splits into isomorphisms between splittings of its domain and codomain, it is an isomorphism. 
\end{proof}
We now define the vector field $f: \B \to \T \B$ by
\begin{equation}\label{eq:fh_def}
f_h(x):= T_x^{-1} \left[T_{P(x)} \D P_x  T_x^{-1}\right]^\dagger T_{P(x)}g(P(x))
\end{equation}

\begin{Prop}\label{prop:specific}
	$f_h$ is a lift of $g$, and $f_h|_\M = g$.
\end{Prop}
\begin{proof}
	We first prove that $f_h$ is a lift. 
	That is, we show that $\forall x \in \B: \D P_x f_h(x) = g(P(x))$.
	The proof is a computation:
	\begin{align*}
	\D P_x f_h(x) &= \D P_x T_x^{-1} \left[T_{P(x)}\D P_x T^{-1}_x\right]^\dagger  T_{P(x)}  g (P(x))\\
	&= T_{P(x)}^{-1}[T_{P(x)} \D P_x T^{-1}_x]T_x \ldots\\
	& \ldots  T^{-1}_x \left[T_{P(x)} \D P_x  T^{-1}_x\right]^\dagger  T _{P(x)}  g (P(x))\\
	&= g(P(x)),
	\end{align*}
	with the last equality following since $\left[T_{P(x)}\D P_x T^{-1}_x\right] \left[T_{P(x)}\D P_x T^{-1}_x\right]^\dagger $ is the identity because $\left[T_{P(x)}\D P_x T^{-1}_x \right]$ is surjective. 
	This proves $f_h$ is a lift of $g$.
	
	We now show that $f_h |_\M = g$.
	Since $P|_\M$ is the identity, we need to show that $\forall x \in \M$:
	\begin{equation*}
	T^{-1}_x  \left[T_x \D P_x T^{-1}_x\right]^\dagger  T_x g(x) = g(x).
	\end{equation*}
	But for every $x \in \M$, $g(x)$ is tangent to $\M$ by assumption; that is, $g(x) \in \ker \D G_x$.
	For $x \in \M$, $\D P_x$ is the projection $\ker \D P_x \oplus \ker \D G_x \to \ker \D G_x = \textnormal{im } \D P_x$.
	Lemma \ref{lem:T_lem} then implies that $\left[T_x \D P_x T^{-1}_x\right]$ is an orthogonal projection map.
	But the Moore-Penrose pseudoinverse of an orthogonal projection is the orthogonal projection itself\footnote{In coordinates, an orthogonal projection may be written as $U \Sigma U^T$ (where $U$ is a square orthogonal matrix and $\Sigma$ is diagonal with diagonal entries equal to $1$ or $0$, and with rank equal to the dimension of the subspace projected onto), so $(U \Sigma U^T)^\dagger = U \Sigma U^T$ using the singular value decomposition formula for the pseudoinverse.}, and hence 
	\begin{align*}
	T^{-1}_x  \left[T_x \D P_x T^{-1}_x\right]^\dagger  T_x (g(x)) 
	&= T^{-1}_x  \left[T_x \D P_x T^{-1}_x\right]  T_x (g(x))
	\\ &= \D P_x g(x)\\
	&= g(x),
	\end{align*}
	since the fact that $P|_\M$ is the identity implies $\D P_x$ is the identity on $\T_x \M$.
\end{proof}

\subsection{Stabilizing $\M$}\label{sec:stabilizing}

We define the $\C^{r-1}$ vector field $f_v: \B \to \T \B$ by:
\begin{equation}\label{eq:fv_def}
f_v(x) = -\alpha(x)\Pi^P_x \nabla V(x), 
\end{equation}
where $\Pi^P$ is defined in equation \eqref{eq:Pi_P_def}, $\alpha_0>0$, and $\alpha: \B \to \R$ is a $\C^{r-1}$ function, such that $\forall x\in\B:\,\alpha(x)\geq\alpha_0$.

\begin{Prop}\label{prop:fv}
	The flow of $f_v$ is complete.
	$\M$ is asymptotically stable under the flow of $f_v$ with basin of attraction $\B$, and each fiber of $P$ is an invariant manifold of $f_v$.
\end{Prop}
\begin{proof}
	First, we show that fibers of $P$ are invariant under the flow.
	By construction, for any $x \in \B$ we see that $f_v$ lies in the tangent space $\T_x P^{-1}(P(x))$ (because $\Pi^P_x$ is a projection onto $\ker \D P_x = \T_x P^{-1}(P(x))$.
	From this it follows that the fibers of $P$ are invariant manifolds of $f_v$.
		
	Our definition of $V$ implies that $\forall x \in \M: \nabla V(x) =0$. Hence $\forall x \in \M: f_v(x) = 0$, making $\M$ an invariant manifold of $f_v$.
	We compute the Lie derivative of $V$ along $f_v$ as follows (using the standard inner product on $\R^n$) for $x \in \B \ct \M$:
	\begin{align*}
	L_{f_v}V(x) &= \langle \nabla V(x),f_h(x)\rangle \\
	&= - \alpha(x) \left \langle \nabla V(x), \Pi^P_x \nabla V(x) \right \rangle.
	\end{align*}
	Since $\Pi^P_x = [I - \D P_x^\dagger \D P_x]$ is a projection operator and projections are always positive-semidefinite, the assumptions on $V$ imply that the right side is zero on $\M$ and strictly negative off of $\M$. 
	It now follows that the flow of $f_v$ is complete, since the trajectory with initial condition $x$ is confined in positive time to the set $V^{-1}(-\infty,V(x)]$ which is compact by Lemma \ref{lem:compact_sublevel}.
	
	Since $\M$ is compact, $V$ is zero on $\M$ and positive on $\B \ct \M$, $V$ does not attain its supremum on $\B$, and (by the assumption in \S \ref{sec:completeness}) $V(x) \to \sup_{y\in \B}V(y)$ as $x$ tends to any point of $\partial \B$, the Lyapunov theorem (\citet{wilson1967structure} Theorem 3.1) implies that $\M$ is asymptotically stable with basin of attraction $\B$.
\end{proof}

\subsection{Properties of the resulting vector field}\label{sec:prop_resulting_f}

As mentioned earlier, we define $$f:= f_v + f_h.$$
In this section, we prove that $f$ has the desired properties.

\begin{Prop}\label{prop:f_props}
$f:\B \to \T \B$ is complete, $\M$ is an asymptotically stable invariant manifold of $f$ with basin of attraction equal to $\B$, and the fibers of $P$ are invariant under the flow of $f$.
\end{Prop}
\begin{proof}
	$f_v(x)$ vanishes on $\M$, so for all $x \in \M$ we see that $f(x) = f_h(x) = g(x)$ is in $ \T _x \M = \ker \D G_x$.
	It follows that $\M$ is an invariant manifold of $f$.
	Next, a computation shows:
	\begin{align*}
	L_f V(x) &= L_{f_v}V(x) + L_{f_h}V(x)\\
	&= L_{f_v}V(x).
	\end{align*}
	$L_{f_h}V$ is zero since $\nabla V(x) = 2 (\D G_x)^TG(x)$\footnote{where $(\cdot)^T$ is the adjoint or transpose} and $f_h(x) \in \Hor_x \B = \ker \D G_x$ imply that $$L_{f_h}V(x) = \langle V(x), f_h(x)\rangle = \langle 2 (\D G_x)^TG(x), f_h(x)\rangle = \langle G(x), \D G_x f_h(x)\rangle = \langle G(x), 0 \rangle = 0.$$
	
    Let $x_0 \in \B$ and let $\phi_t(x_0)$ be the solution at time $t$ of the initial value problem $\dot x = f(x), x(0)=x_0$.
	Since $L_f V = L_{f_v}V < 0$,
	for any $t > 0$ with $t$ in the maximal interval of existence of $\phi_t(x_0)$, $\phi_t(x_0)\in V^{-1}(-\infty,V(x)]$.
	The assumption in \S \ref{sec:completeness} shows that $V^{-1}(-\infty,V(x)]$ is compact, so $\phi_t(x_0)$ is defined for all $t > 0$ and thus $f$ is complete. 
    
    As in the proof of Proposition \ref{prop:fv}: since $\M$ is compact, $V$ is zero on $\M$ and positive on $\B \ct \M$, $V$ does not attain its supremum on $\B$, and (by the assumption in \S \ref{sec:completeness}) $V(x) \to \sup_{y\in \B}V(y)$ as $x$ tends to any point of $\partial \B$, the Lyapunov theorem (\citet{wilson1967structure} Theorem 3.1) implies that $\M$ is asymptotically stable with basin of attraction $\B$.
    
	Finally, for every $x \in \B$ we have $\D P_x f(x) = \D P_x f_h(x) = g(P(x))$, using the result of Proposition \ref{prop:specific} and the fact that $f_v(x) \in \ker \D P_x$.
    It follows that the fibers of $P$ are invariant under the flow of $f$.
	This completes the proof.
\end{proof}

We now set out to prove that $\M$ is (locally) exponentially stable, and furthermore that the rate of exponential convergence can be made arbitrarily large by choosing $\alpha_0$ appropriately.

\begin{Lem}\label{lem:exp_stab_assump}
Given $x\in \B$, let $d_x:=\inf_{y\in\M}\|x-y\|$ denote the distance from $x$ to $\M$.
There exists an open neighborhood $U_E$ of $\M$ and positive constants $k_1, k_2, k_3 > 0$ on which the following holds:

\begin{enumerate}
	\item $\forall x \in U_E: k_1 d_x^2 \leq V(x) \leq k_2 d_x^2$
	\item $\forall x \in U_E: \langle \nabla V(x), \Pi^P_x \nabla V(x)\rangle \geq k_3 d_x^2 $
\end{enumerate}
	
\end{Lem}
\begin{proof}
	$V(x) = \|G(x)\|^2$ by definition, so for $x \in \B$ and $y \in \M$ we have 
	\begin{align*}
	V(x) &= \left[\D G_y(x-y)+\bo(\|x-y\|^2)\right]^2\\
	&= \left \langle x-y,\D G_y^T\D G_y (x-y)\right\rangle + \bo(\|x-y\|^3),
	\end{align*}
	since $G$ vanishes on $\M$.
    Here, $\bo(\cdot)$ is the standard \concept{big-oh} notation.
	For any $y\in \M$, $\D G_y^T \D G_y$ maps $\Nor_y\M$ isomorphically onto itself.\footnote{To see this, note that $\ker \D G_y = \T_y\M$, so $\T_y\M \subseteq \ker \D G_y^T \D G_y$. 
	To show the reverse inclusion, if $\D G_y^T\D G_y v = 0$, then $\langle \D G_y^T\D G_y v,v\rangle = \langle \D G_y v, \D G_y v\rangle = 0$, so $v \in \ker \D G_y = \T_y \M$.
	Hence the rank of $\D G_y$ is equal to the dimension of $\Nor_y\M$.
	Finally, the range of $\D G_y^T \D G_y$ is contained in $\Nor_y\M$, since for any $v \in \T_y\R^n, w \in \T_y\M$, we have 
	$\langle \D G_y^T\D G_y v,w\rangle = \langle \D G_y v,\D G_y w\rangle = \langle \D G_y v,0\rangle = 0$, so the rank-nullity theorem implies that the range of $\D G_y^T \D G_y$ is actually all of $\Nor_y\M$ and hence for $y \in \M$, $\D G_y^T \D G_y$ maps $\Nor_y\M$ isomorphically onto itself. }
    From this and compactness of $\M$ it follows that $a:= \min_{y\in\M}\|\D G_y^T \D G_y|_{\Nor_y\M}\| > 0$.
    We also have $A:=\max_{y\in\M}\|\D G_y^T \D G_y\| < \infty$.
    Using these facts and choosing $y\in\M$ to minimize the distance from $x$ to $\M$, we find (since $(x-y)\in\Nor_y\M$ in this case):
    \begin{equation*}
    a\|x-y\|^2 - \bo(\|x-y\|^3) \leq V(x) \leq A\|x-y\|^2 + \bo(\|x-y\|^3),
    \end{equation*}	
    or 
    \begin{equation*}
    a d_x^2 - \bo(d_x^3) \leq V(x) \leq A d_x^2 + \bo(d_x^3). 
    \end{equation*}
Choose $U_1$ so that $\sup_{x\in U_1}d_x$ is sufficiently small so that $|\bo(d_x^3)| < \frac{a}{2}d_x^2$.
    Taking $k_1:= a-\frac{a}{2}$ and $k_2:= A + \frac{a}{2}$ completes the first part of the proof.
    
    Since $\Pi^P_x\nabla V(x)\not = 0$ for all $ x \in \B$, taking $U_2$ to be a precompact neighborhood of $\M$
    shows that $$b:= \inf_{x\in U_2}\|\Pi^P_x|_{\text{span}\{\nabla V(x)\}}\|> 0.$$
    It follows that 
    \begin{align*}
    \langle \nabla V(x), \Pi^P_x \nabla V(x)\rangle &\geq b \|\nabla V(x)\|^2.
    \end{align*}
    For $x \in \B$ and $y \in \M$ we have
    \begin{equation*}
    V(x) = \nabla V(x)(x-y)+\bo(\|x-y\|^2).
    \end{equation*}
    Squaring this equation yields
    \begin{equation*}
    k_1^2 d_x^4 \leq V^2(x) = \langle \nabla V(x)(x-y), \nabla V(x)(x-y)\rangle + \bo(\|x-y\|^3), 
    \end{equation*}
    where the first inequality follows from the proof of the first claim.
    Choosing $y\in \M$ so that $\|x-y\| = d_x$ and taking norms, we have 
    \begin{equation*}
     k_1^2d_x^4 \leq \|\nabla V(x)\|^2d_x^2 + \bo(d_x^3)
    \end{equation*}
    Choose $U_2$ to be sufficiently small that $\sup_{x\in U_2}d_x < 1$ and $|\bo(d_x^3)| < \frac{k_1^2}{2}d_x^2$ for all $x \in U_2$.
    It follows that
    \begin{equation*}
    \langle \nabla V(x), \Pi^P_x \nabla V(x)\rangle \geq b \|V(x)\|^2 \geq  k_3 d_x^2
    \end{equation*}
    for $x \in U_2$, with $k_3:=b\left(\frac{k_1^2}{2}\right)$.
    Taking $U_E:=U_1\cap U_2$ completes the proof.
\end{proof}

\begin{Prop}\label{prop:exp_stable}
$\M$ is exponentially stable on an open neighborhood $U_E$ of $\M$.
	Furthermore, the rate of exponential convergence can be made arbitrarily large by picking $\inf_{x\in \bar U_E} \alpha(x)$ large.
\end{Prop}
\begin{proof} 

Let $U_E,k_1,k_2,k_3,a$ be as in Lemma \ref{lem:exp_stab_assump} and define $k_4:= \inf_{x \in \bar U_E}\alpha(x)$.
	For all $x \in U_E$, the conclusion of Lemma \ref{lem:exp_stab_assump} implies that
	
	\begin{align*}
	\frac{d}{dt}V(\phi_t(x)) &= \langle \nabla V(\phi_t(x)), f_v(\phi_t(x))+ f_h(\phi_t(x))\rangle \\
    &= \langle \nabla V(\phi_t(x)), -\alpha(\phi_t(x))\Pi^P_{\phi_t(x)} \nabla V(\phi_t(x)) + f_h(\phi_t(x))\rangle \\
    &= \langle \nabla V(\phi_t(x)), -\alpha(\phi_t(x))\Pi^P_{\phi_t(x)} \nabla V(\phi_t(x))\rangle\\ 
    &\leq - \alpha(\phi_t(x))\frac{k_3}{k_2}V(\phi_t(x)) \\
    &\leq -\frac{k_4 k_3}{k_2}V(\phi_t(x)),
	\end{align*}
	with the third equality following since $\langle \nabla V(\phi_t(x)), f_h(\phi_t(x)\rangle = 0$ as explained in the proof of Proposition \ref{prop:f_props}.
	Again letting $d_x$ denote the distance from $x$ to $\M$,  the result of Lemma \ref{lem:exp_stab_assump} together with an application of Gronwall's inequality together yields
	
	\begin{align*}
	d_{\phi_t(x)} & \leq \left[ \frac{V(\phi_t(x))}{k_1}\right]^{\frac{1}{2}} \\
	&\leq \left[\frac{V(x)}{k_1}\right]^{\frac{1}{2}}e^{-\frac{k_4 k_3}{2 k_2}t}\\
	& \leq C d_x e^{-\mu t},
	\end{align*}
	with $C = \left[\frac{k_2}{k_1}\right]^\frac{1}{2} > 0$ and $\mu = \frac{k_4 k_3}{2 k_2}.$
	We see that we can make $\mu$ arbitrarily large by picking $k_4 = \inf_{x \in \bar U_E}\alpha(x)$ sufficiently large.
    This completes the proof.
\end{proof}

\section{Normal hyperbolicity of $\M$ and robustness of our construction}\label{sec:our_persistence}

We have shown that under the flow induced by the vector field $f$ on $\B$, $\M$ is asymptotically stable with basin of attraction equal to $\B$.
We have also shown that $\M$ is exponentially stable on the neighborhood $U_E \supset \M$, with exponential rate $\mu$ proportional to $\min_{x \in \bar U_E}\alpha(x)$. 
In Appendix \ref{app:proofs_persistence}, we show that if $\min_{x \in \bar U_E}\alpha(x)$ (and hence $\mu$) is chosen sufficiently large, $\M$ can be made $k$-normally hyperbolic for any $k \in \N$.
%\SR{is this true? or is $k$ limited by the dimension of the space? I'm confused...}
%\MK{This is true. $k$ is the factor by which normal contraction dominates tangent contraction, to first order. It has nothing to do with the dimension of the space. In Fenichel's version of NHIM theory $k$ can be any real number, even a negative real. However, if the NHIM $\M$ is a $\C^k$ manifold, nothing useful really happens if $M$ is $r$-normally hyperbolic with $r > k$; $\M$ still just persists to a $\C^k$ manifold, etc.}
As a corollary of this fact and other results in Appendix \ref{app:proofs_persistence}, we have the following Theorem \ref{th:persistence_of_our} showing that our construction in \S \ref{sec:construction} is robust -- it persists under perturbations.
We prove this theorem in Appendix \ref{app:proofs_persistence}.

\begin{Th}\label{th:persistence_of_our}
	Assume $r > 3$.
	Let $\mu = \frac{k_4 k_3}{2 k_2}$ be as in Proposition \ref{prop:exp_stable}, where $k_4 = \min_{x \in \bar U_E}\alpha(x)$,
	and choose $\alpha:\B\to \R$ so that 
	$$k_4 > r \frac{2 k_2}{k_1}L.$$
	Then there exists $\theta > 0$ sufficiently small such that if $g:\B \to \T\B$ is another $\C^{r-1}$ vector field such that
	\begin{align*}
	\sup_{x\in\B}\|g(x)-f(x)\| &< \theta\\
	\sup_{x\in\B}\|\D g(x)-\D f(x)\|  &< \theta,
	\end{align*}	
    then there exists an open set $\B^g\subseteq \B$ positively invariant under the flow of $g$ and a $\C^{r-1}$ exponentially stable normally hyperbolic submanifold $\M^g$ $\C^{r-1}$ diffeomorphic to $\M$ and $\C^1$-close to $\M$.
    The stability basin of $\M^g$ contains $\B^g$.
    $\M^g$ has the unique asymptotic phase property with a $\C^{r-2}$ phase map $P^g:\B^g\to\M^g$ making $(\B^g,P^g,\M^g)$ into a $\C^{r-2}$ fibered manifold with $(n-k)$-dimensional Euclidean fibers.
    The fibers of $P^g$ are $\C^1$-close to the fibers of $P$ on $\B^g$.	
\end{Th}

\section{Examples}\label{sec:examples}
\subsection{Basic examples}
In this section we present some simple examples which we hope nonetheless serve as concrete illustrations of the theory.
In all examples, we created plots of trajectories using a \verb|NumPy| implementation \citep{integro_py} of the dopri5 ODE integrator \citep{hairer2010solving} to numerically integrate vector fields.
\subsubsection{An equilibrium point}
Let $\Q = \B = \R^2$ with the Euclidean inner product and let $\M = \{\{0,0\}\}$ and define $g:\M\to\T \M$ by $g(0,0)=[0,0]^T$.
$\M$ is the zero level set of the smooth submersion $G:\R^2 \to \R^2$ defined by $G(x,y)=[x,y]^T$.
The function $P:\R^2 \to \R^2$ is defined by $P(x,y)=[0,0]^T$ for all $(x,y)\in\R^2$.
It follows that $\forall (x,y)\in \R^2$, $\ker \D P_{(x,y)} = \T_{(x,y)}\R^2$ and $\ker \D G_{(x,y)} = \{(0,0)\}\in \T_{(x,y)}\R^2$.
Hence $\R^2 = \ker \D G \oplus \ker \D P$, so the assumption in \S \ref{sec:connection} is satisfied.
Additionally, $\|G\|$ tends to $\infty$ as $\|(x,y)\| \to \infty$, so the assumption in \S \ref{sec:completeness} is satisfied.
$f_h:\B\to \T \B$ as defined in equation \eqref{eq:fh_def} is the zero vector field, and we thus have $f = f_v$.
$\Pi^P_{(x,y)}$ in the definition of $f_v$ in equation \eqref{eq:fv_def} is the identity map on $\T_{(x,y)}\R^2$ and thus (taking the function $\alpha$ to be $\alpha \equiv 1/2$) $f = f_v$ is given by
\begin{align*}
f(x,y) &= - \frac{1}{2}\nabla \|G(x,y)\|^2\\
&= -[x,y]^T.
\end{align*}

\subsubsection{A limit cycle}\label{sec:ex-lim-cyc}
Let $\Q = \R^2$ with the Euclidean inner product, let $\B = \{(x,y)\in\R^2| \|(x,y)\|\leq 2\} \ct \{(0,0)\}$, and let $\M = S^1 := \{(x,y)\in \R^2|x^2+y^2 = 1\}$.
Note that $\M$ is the zero level set of the smooth submersion $G:\B \to \R$ defined by $G(x,y):= x^2 + y^2 - 1$, and $\|G\|$ tends to the constant value $1$ as $(x,y) \to \partial \B$ from within $\B$ in accordance with the completeness assumption\footnote{In this example, everything would actually work fine if we took $\B:= \R^2 \ct \{0,0\}$. However, our assumption in \ref{sec:completeness} does not guarantee that this would be the case because then $\|G\|$ would not approach a constant value on $\partial \B \cup \{\infty\}$. Perhaps a better assumption in lieu of the one in \S \ref{sec:completeness} would eliminate this technical annoyance.} in \S \ref{sec:completeness}.
Define $P:\B \to \B$ by $P(x,y):= \frac{[x,y]^T}{\sqrt{x^2 + y^2}}$.
Define $g:\M \to \T \M$ by $g(x,y):= [-y,x]^T$.
We compute 
%\iffalse 
\begin{equation*}
\D P_{(x,y)} = \frac{1}{\sqrt{x^2+y^2}}
\begin{bmatrix}
1-\frac{x^2}{x^2+y^2} & -\frac{xy}{x^2+y^2}\\
-\frac{xy}{x^2+y^2} & 1-\frac{y^2}{x^2+y^2}
\end{bmatrix}
\end{equation*}
%\fi
By inspection, $\ker \D P_{(x,y)} = \text{span}\{[x,y]^T\}$ and $\ker \D G_{(x,y)} = \text{span}\{[x,y]^T\}^\perp = \text{span}\{[-y,x]^T\}$, so $\T\B = \ker \D G \oplus \ker \D P$ in accordance with the assumption in \S \ref{sec:connection}.
Using the requirements that $\D P_{(x,y)}f_h(x,y) = g(P(x,y))$ and $f_h(x,y)\in \ker \D G_{(x,y)}$, we find:

\begin{equation*}
f_h(x,y) = \left[-y,x\right]^T,
\end{equation*}
since $\D P_{(x,y)}[-y,x]^T = \frac{1}{\sqrt{x^2+y^2}}[-y,x]^T = g\left(\frac{x}{\sqrt{x^2+y^2}},\frac{y}{\sqrt{x^2+y^2}}\right) = g(P(x,y))$.

By the chain rule, $\nabla V(x,y) = 2\D G^T_{(x,y)}G(x,y) = 2(x^2+y^2-1)[2x,2y]^T$ and $\Pi^P_{(x,y)}$ is the identity when restricted to the subspace $\text{span}\{[x,y]^T\}$. 
We thus compute $f_v(x,y):= -\alpha(x,y) \Pi^P_{(x,y)} \nabla V(x,y)$ as:

\begin{equation*}
f_v(x,y) = -4(x^2+y^2-1)\alpha(x,y)\left[x,y\right]^T.
\end{equation*}

\begin{figure}[h]
	\centering 
	\includegraphics[width=1.0\textwidth]{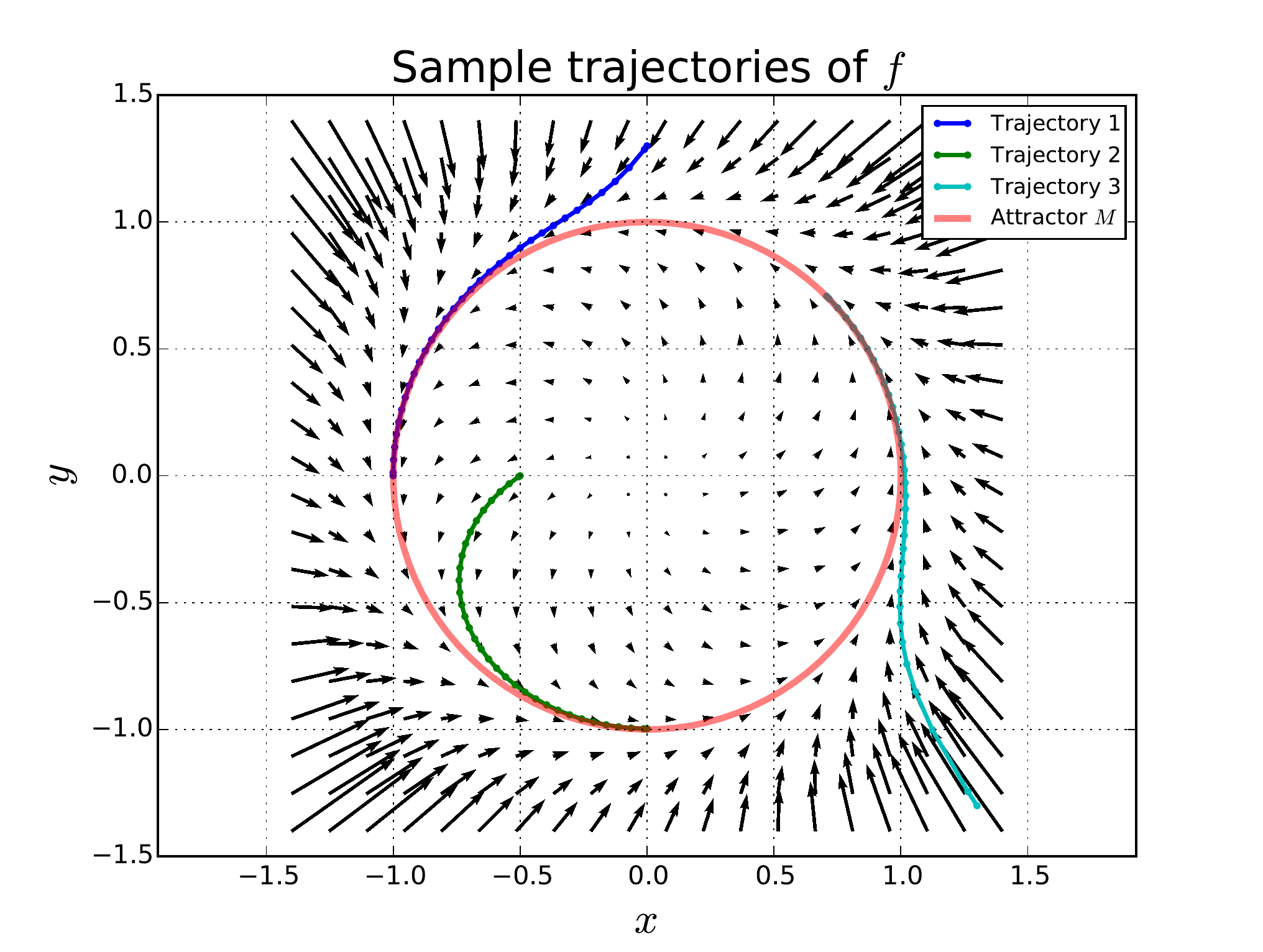}
	\caption[]{A quiver plot of $f:= f_h + f_v$ from the example in \S \ref{sec:ex-lim-cyc} is shown together with three sample trajectories of $f$.
	The attractor, $M$, is shown in red.}
	\label{fig:ex-lim-cyc}
\end{figure}

We now form the vector field $f:= f_v + f_h$.
We chose $\alpha(x,y)\equiv 0.5$ and plotted the resulting the vector field $f$ and multiple trajectories of $f$ in Figure \ref{fig:ex-lim-cyc}.

\begin{figure}[h]
	\centering 
	\includegraphics[width=1.0\textwidth]{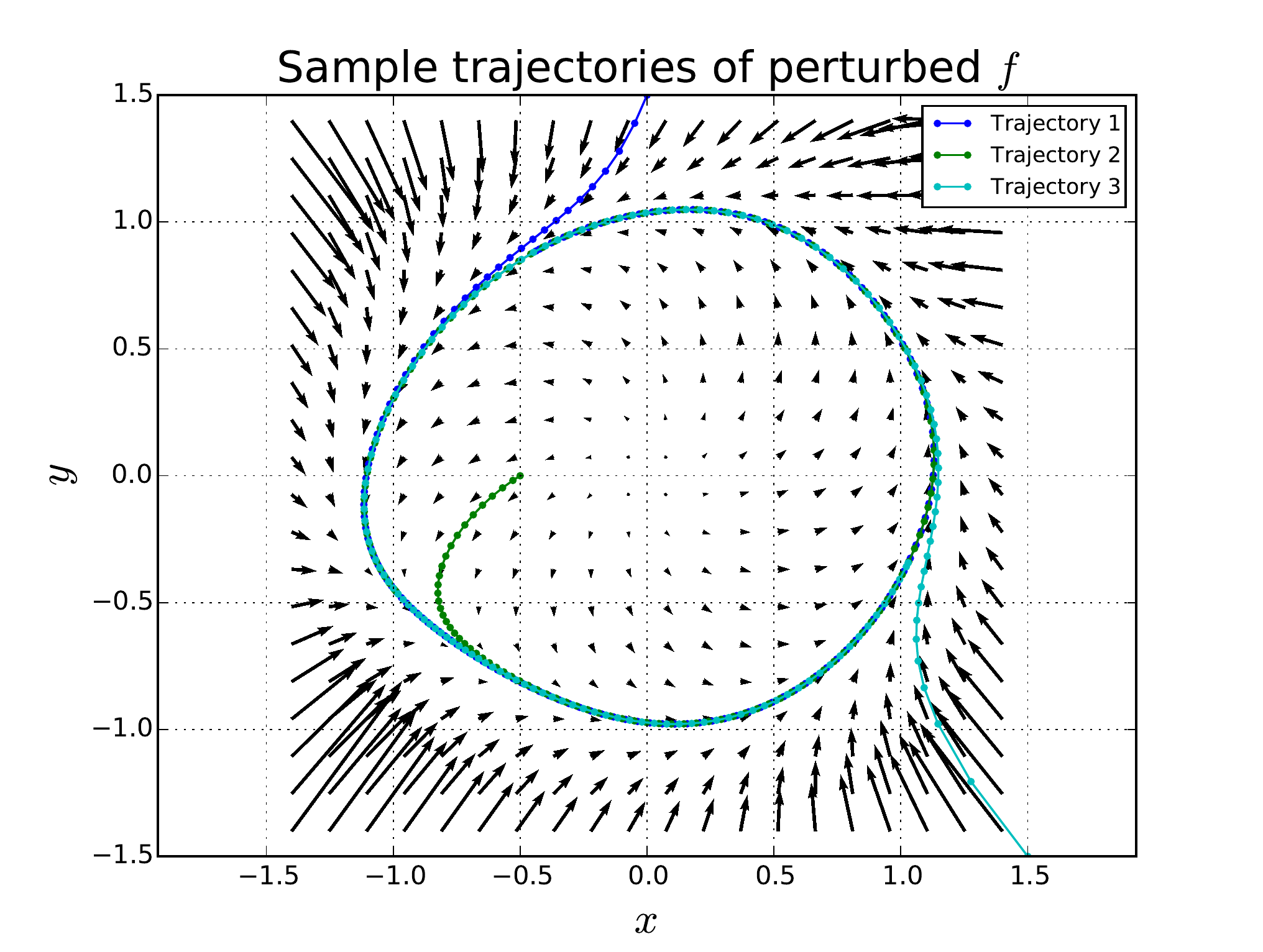}
	\caption[]{A quiver plot of the perturbed $f$ from the example in \S \ref{sec:ex-lim-cyc} is shown together with three sample trajectories of this perturbed vector field.
    This plot suggests that the limit cycle $\M$ persists under the perturbation to a new limit cycle $\tilde \M$, consistent with the conclusion of our Theorem \ref{th:persistence_of_our}.}
	\label{fig:ex-lim-cyc-pert}
\end{figure}

Next, we illustrate the robustness of the structure of our construction to perturbations of the vector field $f$.
We form the perturbed vector field $f^{\text{pert}}$ as follows:

\begin{equation*}
f^{\text{pert}}(x,y):= f(x,y) + \varepsilon [\eta_1(x,y),\eta_2(x,y)]^T,
\end{equation*}
where $\eta_1,\eta_2: \B \to \R$ are $\C^2$ functions and $\varepsilon > 0$ is a small parameter.
Theorem \ref{th:persistence_of_our} says that for $\varepsilon > 0$ sufficiently small, $\M$ and $P$ persist -- $\M$ is deformed into a $\C^1$-close invariant manifold $\tilde \M$ diffeomorphic to $\M$, and $P$ is deformed into a $\C^1$-close phase map $\tilde P$.
We arbitrarily chose to define $\eta_1(x,y):= x^3\cos(xy)$ and $\eta_2(x,y):=xye^{xy}$.
Trajectories of $f^{\text{pert}}$ are shown in Figure \ref{fig:ex-lim-cyc-pert} for $\varepsilon = 0.5$, which illustrates the persistence of $\M$.

\subsubsection{An invariant sphere}\label{sec:ex-sphere}
Let $\Q = \R^3$ with the Euclidean inner product, let $\B = \{(x,y,z)|\|(x,y,z)\leq 2\} \ct \{(0,0,0)\}$, and let $\M = S^2 := \{(x,y,z)\in \R^3|x^2+y^2+z^2 = 1\}$.
Note that $\M$ is the zero level set of the smooth submersion $G:\B \to \R$ defined by $G(x,y,z):= x^2 + y^2 +z^2 - 1$, and $\|G\|$ tends to the constant value $1$ as $(x,y,z) \to \partial \B$ from within $\B$ in accordance with the completeness assumption\footnote{In this example (similarly to the last example), everything would actually work fine if we took $\B:= \R^3 \ct \{0,0\}$. However, our assumption in \ref{sec:completeness} does not guarantee that this would be the case because then $\|G\|$ would not approach a constant value on $\partial \B \cup \{\infty\}$. Perhaps a better assumption in lieu of the one in \S \ref{sec:completeness} would eliminate this technical annoyance.} in \S \ref{sec:completeness}.
Define $P:\B \to \B$ by $P(x,y,z):= \frac{[x,y,z]^T}{\sqrt{x^2 + y^2 + z^2}}$.
Define $g:\M \to \T \M$ by $g(x,y,z):= [-y,x,0]^T$.
The choice of $\M, P$, and $G$ will make the analysis very similar to that of the example in \S \ref{sec:ex-lim-cyc}. 
We compute 

\begin{equation*}
\D P_{(x,y,z)} = \frac{1}{\sqrt{x^2+y^2+z^2}}
\begin{bmatrix}
1-\frac{x^2}{x^2+y^2+z^2} & -\frac{xy}{x^2+y^2+z^2} & -\frac{xz}{x^2+y^2+z^2}\\
-\frac{xy}{x^2+y^2+z^2} & 1-\frac{y^2}{x^2+y^2+z^2} & -\frac{yz}{x^2+y^2+z^2}\\
-\frac{xz}{x^2+y^2+z^2} & -\frac{yz}{x^2+y^2+z^2} & 1-\frac{z^2}{x^2+y^2+z^2}
\end{bmatrix}.
\end{equation*}
By inspection, $\ker \D P_{(x,y,z)} = \text{span}\{[x,y,z]^T\}$ and $\ker \D G_{(x,y,z)} = \text{span}\{[x,y,z]^T\}^\perp$, so $\B = \ker \D G \oplus \ker \D P$ in accordance with the assumption in \S \ref{sec:connection}.
Using the requirements that $\D P_{(x,y,z)}f_h(x,y,z) = g(P(x,y,z))$ and $f_h(x,y,z)\in \ker \D G_{(x,y,z)}$, we find:

\begin{equation*}
f_h(x,y) = \left[-y,x,0\right]^T,
\end{equation*}
since $\D P_{(x,y)}[-y,x,0]^T = \frac{1}{\sqrt{x^2+y^2+z^2}}[-y,x,0]^T = g\left(\frac{x}{\sqrt{x^2+y^2+z^2}},\frac{y}{\sqrt{x^2+y^2+z^2}},\frac{z}{\sqrt{x^2+y^2+z^2}}\right) = g(P(x,y,z))$.

By the chain rule, $\nabla V(x,y,z) = 2\D G^T_{(x,y,z)}G(x,y,z) = 2(x^2+y^2+z^2-1)[2x,2y,2z]^T$ and $\Pi^P_{(x,y,z)}$ is the identity for all $(x,y,z)$. 
We thus compute $f_v(x,y,z):= -\alpha(x,y,z) \Pi^P_{(x,y,z)} \nabla V(x,y,z)$ as:

\begin{equation*}
f_v(x,y,z) = -4(x^2+y^2+z^2-1)\alpha(x,y,z)\left[x,y,z\right]^T.
\end{equation*}

\begin{figure}[h]
	\centering 
	\includegraphics[width=1.0\textwidth]{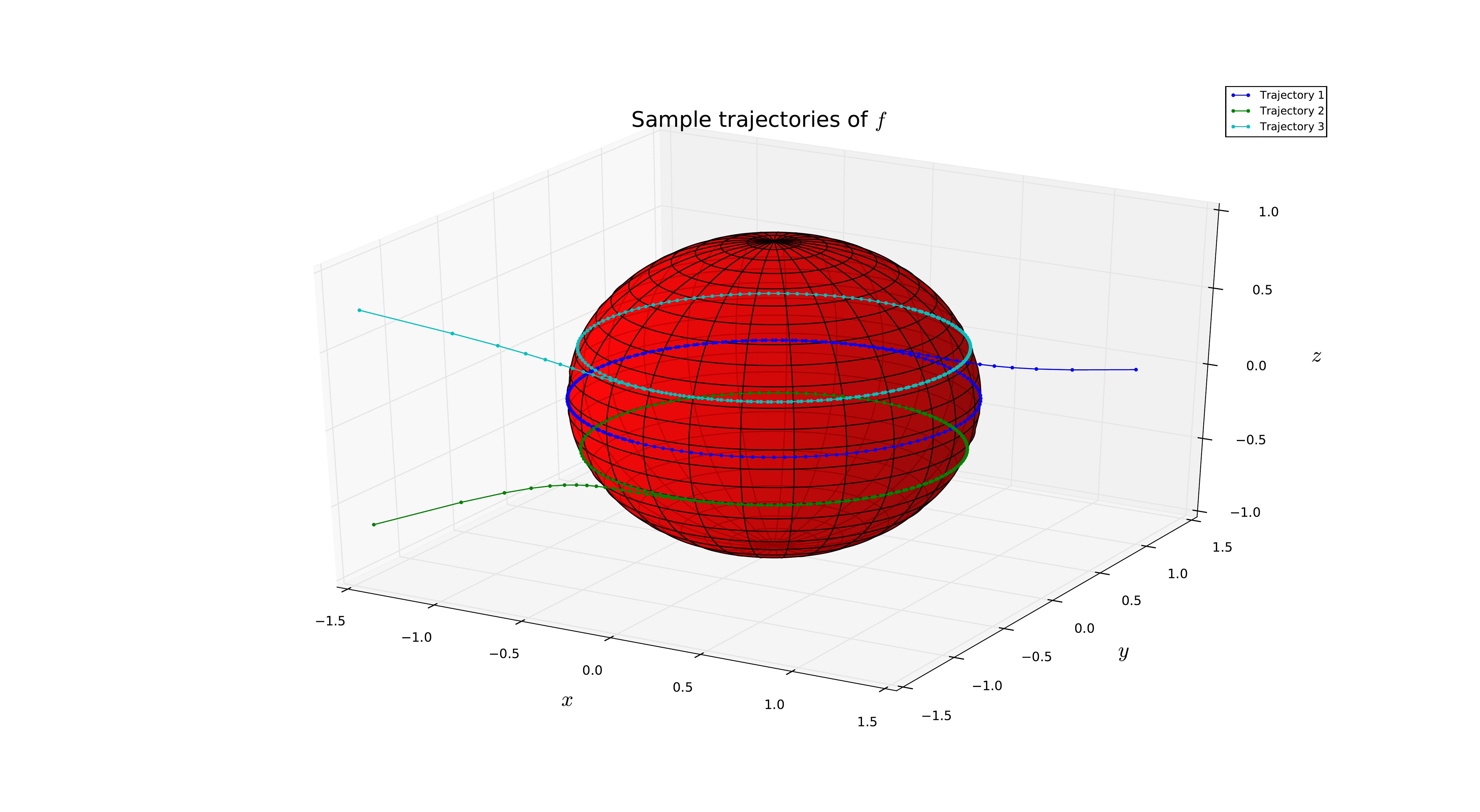}
	\caption[]{Three sample trajectories of $f$ from the example in \S \ref{sec:ex-sphere} are shown.
	The attractor, $M$, is shown in red.}
	\label{fig:ex-sphere}
\end{figure}

We now form the vector field $f:= f_v + f_h$.
We chose $\alpha(x,y)\equiv 0.5$ and plotted the resulting vector field $f$ and multiple trajectories of $f$ in Figure \ref{fig:ex-sphere}.

\begin{figure}[h]
	\centering 
	\includegraphics[width=1.0\textwidth]{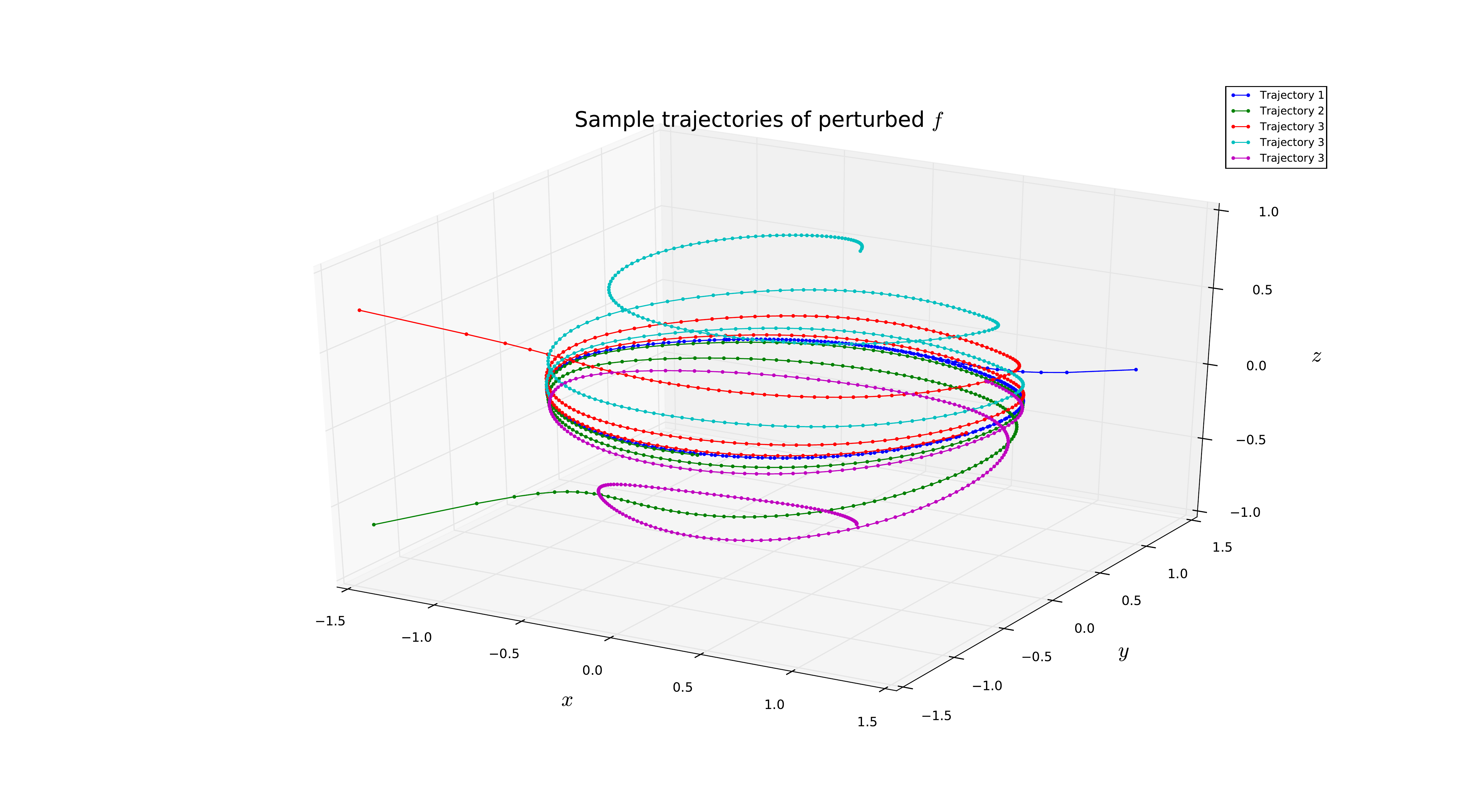}
	\caption[]{Three sample trajectories of the perturbed vector field from from the example in \S \ref{sec:ex-sphere} are shown here.
    This plot suggests that the invariant sphere persists under the perturbation to a new invariant manifold diffeomorphic to a sphere, consistent with the conclusion of our Theorem \ref{th:persistence_of_our}.}
	\label{fig:ex-sphere-pert}
\end{figure}

Next, we illustrate the robustness of the structure of our construction to perturbations of the vector field $f$.
We form the perturbed vector field $f^{\text{pert}}$ as follows:

\begin{equation*}
f^{\text{pert}}(x,y):= f(x,y) + \varepsilon [\eta_1(x,y,z),\eta_2(x,y,z),\eta_3(x,y,z)]^T,
\end{equation*}
where $\eta_1,\eta_2,\eta_3: \B \to \R$ are $\C^2$ functions and $\varepsilon > 0$ is a small parameter.
Theorem \ref{th:persistence_of_our} says that for $\varepsilon > 0$ sufficiently small, $\M$ and $P$ persist -- $\M$ is deformed into a $\C^1$-close invariant manifold $\tilde \M$ diffeomorphic to $\M$, and $P$ is deformed into a $\C^1$-close phase map $\tilde P$.
We arbitrarily chose to define $\eta_1(x,y,z):= x e^{x^2y}\cos(y)$, $\eta_2(x,y,z):=xe^{-z^2}\sin(z)$, and $\eta_3(x,y,z):=xyz$.
Trajectories of $f^{\text{pert}}$ are shown in Figure \ref{fig:ex-sphere-pert} for $\varepsilon = 0.7$, which illustrates the persistence of $\M$.

\subsection{An extended example}
In this section we use the tools of our theory in a more involved example motivated from a physical system -- the double pendulum.
In this extended example, we again created plots of trajectories using a \verb|NumPy| implementation \citep{integro_py} of the dopri5 ODE integrator \citep{hairer2010solving} to numerically integrate vector fields..
\subsubsection{The kinematic double pendulum}\label{sec:kin_pendulum}
Let $\Q := (-\pi,\pi)\times(-\pi,\pi)\subset \R^2$.
Let $\tilde \M = [-\pi,\pi]$ with $\{-\pi\}$ identified with $\{\pi\}$, so that $\tilde \M \cong S^1$.
Note that $\T \tilde \M \cong \tilde \M \times \R$.
Let the vector field $\tilde g:\tilde M \to \T \tilde \M$  be given by $\tilde g(\varphi) = s_h$, with $s_h > 0$ ($h$ for ``horizontal''),
so the dynamics on $\tilde \M$ are given by $\dot \varphi = s_h$.
$\varphi$ should be thought of as the ``phase'' of the oscillation of a single energy-conserving pendulum.
Let $0 < \delta < 2\pi$ and $0 < a_1,a_2 < \pi$ and define the embedding $ F: \tilde \M \to \Q$ by 

\begin{equation*}
F(\varphi) = (\theta_1(\varphi),\theta_2(\varphi)):= (a_1\sin (\varphi-\delta), a_2\sin (\varphi + \delta)),
\end{equation*}
Define $\M = F(\tilde \M) \subseteq \Q$.
Using the identity $\sin(\kappa+\beta)=\sin\kappa\cos\beta+\cos\kappa\sin\beta$, we see that 

\begin{align*}
\frac{\theta_2}{a_2}-\frac{\theta_1}{a_1} & = 2 \sin\delta \cos \varphi\\
\frac{\theta_2}{a_2}+\frac{\theta_1}{a_1} & = 2 \cos \delta \sin \varphi.
\end{align*}
It follows that $\M$ is the zero level set of the smooth map $G:(-\pi,\pi)\times(-\pi,\pi) \to \R$ defined by
\begin{equation}\label{eq:example_G}
 G(\theta_1,\theta_2) = \left[\frac{1}{2\sin \delta}\left( \frac{\theta_2}{a_2}-\frac{\theta_1}{a_1}\right)\right]^2 + \left[\frac{1}{2 \cos \delta}\left(\frac{\theta_1}{a_1} +\frac{\theta_2}{a_2}\right)\right]^2 - 1.
\end{equation}
We compute:
\begin{equation}\label{eq:example_DG}
\D G_{(\theta_1,\theta_2)} =
\begin{bmatrix}
-\frac{1}{2 a_1 \sin^2 \delta}\left(\frac{\theta_2}{a_2}-\frac{\theta_1}{a_1}\right) + \frac{1}{2 a_1 \cos^2 \delta}\left(\frac{\theta_2}{a_2}+\frac{\theta_1}{a_1}\right),  &
\frac{1}{2 a_2 \sin^2 \delta}\left(\frac{\theta_2}{a_2}-\frac{\theta_1}{a_1}\right) + \frac{1}{2 a_2 \cos^2 \delta}\left(\frac{\theta_2}{a_2}+\frac{\theta_1}{a_1}\right) 
\end{bmatrix},
\end{equation}
which is zero if and only if $\frac{1}{\sin^2\delta}\left(\frac{\theta_2}{a_2}-\frac{\theta_1}{a_1}\right) = \frac{1}{\cos^2\delta}\left(\frac{\theta_2}{a_2}+\frac{\theta_1}{a_1}\right)$ and $\frac{1}{\sin^2\delta}\left(\frac{\theta_2}{a_2}-\frac{\theta_1}{a_1}\right) = -\frac{1}{\cos^2\delta}\left(\frac{\theta_2}{a_2}+\frac{\theta_1}{a_1}\right)$, which is possible if and only if $\theta_1 = \theta_2 = 0$.
It follows that $G$ is a submersion on $(-\pi,\pi)\times(-\pi,\pi) \ct \{0,0\}$.
We define $\B$ to be any open neighborhood of $\M$ contained in $(-\pi,\pi)\times (-\pi,\pi)\ct\{0,0\}\subset \Q$ such that $\|G\|$ does not attain its supremum on $\B$ and $\|G\|$ tends to its supremum as $(\theta_1,\theta_2)$ approaches any point of $\partial \B$\footnote{Many such sets $\B$ always exist. We make no effort to explicitly determine a specific $\B$ here.}.
We define the map $\tilde P:\B \to \tilde \M$ simply by extending the formula for $F^{-1}:\M \to \tilde \M$ to all of $\B$.
I.e., $\tilde P$ is given by
\begin{equation}\label{eq:example_tildeP}
\tilde P(\theta_1,\theta_2) = \text{atan2}\left(\frac{1}{2\cos \delta}\left(\frac{\theta_1}{a_1}+\frac{\theta_2}{a_2}\right),\frac{1}{2\sin \delta}\left(\frac{\theta_2}{a_2}-\frac{\theta_1}{a_1}\right)\right).
\end{equation}
For notational brevity, define
\begin{align*}
D&:= \sqrt{\frac{1}{4 \cos^2\delta }\left(\frac{\theta_1}{a_1}+\frac{\theta_2}{a_2}\right)^2 + \frac{1}{4\sin^2 \delta}\left(\frac{\theta_2}{a_2}-\frac{\theta_1}{a_1}\right)^2}\\
k_1&:= \frac{1}{2 a_1 \cos^2\delta }\left(\frac{\theta_1}{a_1}+\frac{\theta_2}{a_2}\right) - \frac{1}{2a_1\sin^2 \delta}\left(\frac{\theta_2}{a_2}-\frac{\theta_1}{a_1}\right)\\ 
k_2&:= \frac{1}{2 a_2 \cos^2\delta }\left(\frac{\theta_1}{a_1}+\frac{\theta_2}{a_2}\right) + \frac{1}{2a_2\sin^2 \delta}\left(\frac{\theta_2}{a_2}-\frac{\theta_1}{a_1}\right).
\end{align*}
Using the identity $\sin(\kappa + \beta) = \sin \kappa \cos \beta + \sin \beta \cos \kappa$, it now follows that $P:\B\to \B$, $P:= F \circ \tilde P$ is given by
\begin{equation*}
P(\theta_1,\theta_2) = \frac{\left(\theta_1, \theta_2\right)}{D},
\end{equation*} 
and
\begin{align*}
\D P_{(\theta_1,\theta_2)} &= 
\frac{1}{D^3}
\begin{bmatrix}
D - \frac{\theta_1}{2 D}k_1 & -\frac{\theta_1}{2 D}k_2\\
-\frac{\theta_2}{2 D}k_1 & D - \frac{\theta_2}{2 D}k_2
\end{bmatrix}.
\end{align*}
We also compute:

\begin{equation*}
\D \tilde P_{(\theta_1,\theta_2)} = \frac{1}{D^2}
\begin{bmatrix}
-\frac{1}{2 \sin \delta}\left(\frac{\theta_2}{a_2} - \frac{\theta_1}{a_1}\right), & \frac{1}{2 \cos \delta}\left(\frac{\theta_2}{a_2} + \frac{\theta_1}{a_1}\right)
\end{bmatrix}.
\end{equation*}

\begin{figure}[h]
	\centering 
	\includegraphics[width=1.0\textwidth]{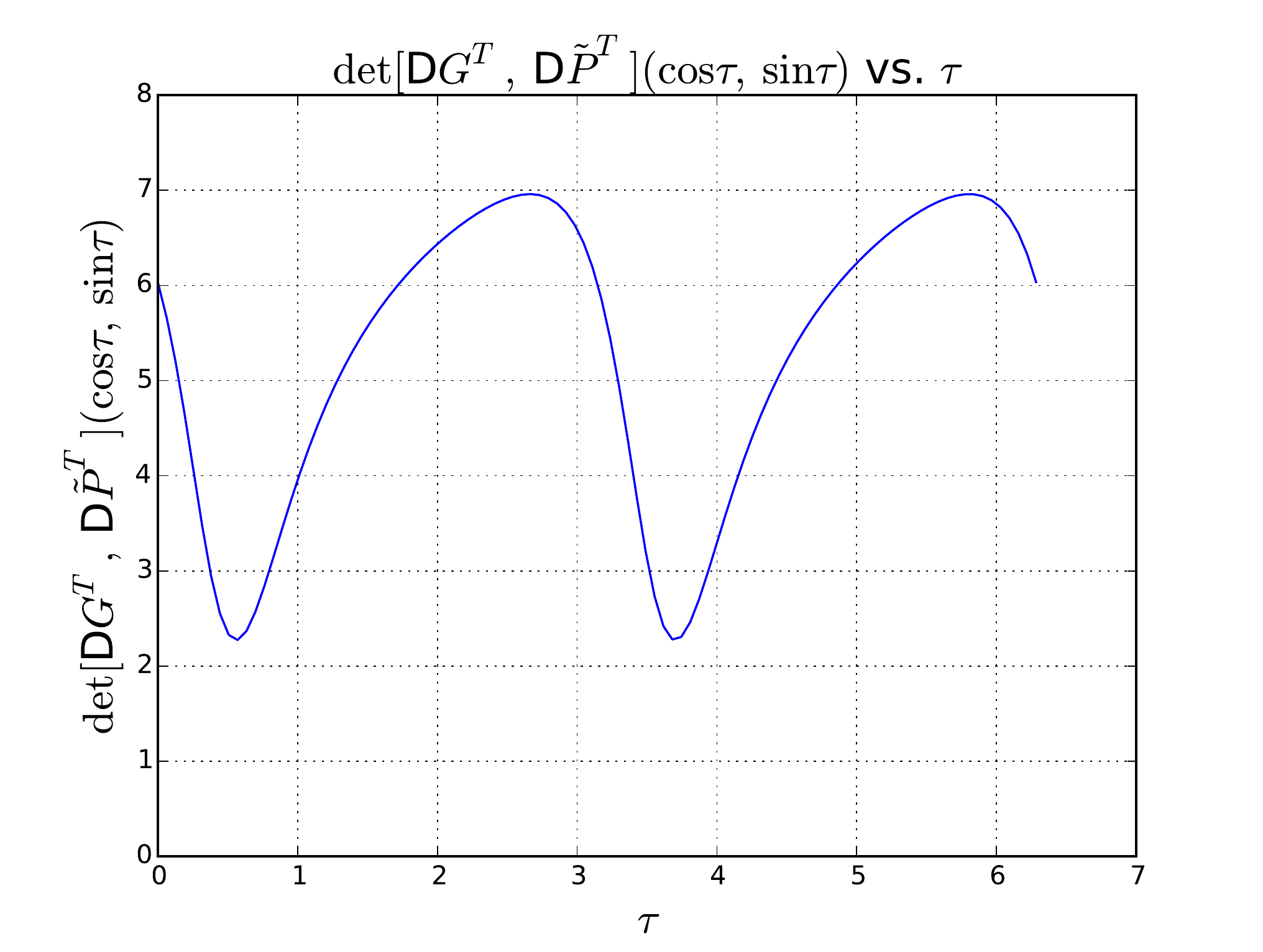}
	\caption[]{A plot of $\det[\D G^T, \D\tilde P^T]$ over the set of points $(\theta_1,\theta_2)$ for which $\theta_1^2 + \theta_2^2 = 1$ is shown.
		Since this determinant is nonzero for any values of $\tau$, it follows that $\T \B = \ker \D G \oplus \ker \D P$.}
	\label{fig:det_plot}
\end{figure}

We now investigate whether the transversality condition \eqref{eq:DG_conn} ($\T \B = \ker \D G \oplus \ker \D P$) holds on $\B$.
Since $F$ is an embedding and $P = F \circ \tilde P$, we see that $\ker \D \tilde P = \ker \D P$.
It follows that for any $(\theta_1,\theta_2)\in \B$, $\T_{(\theta_1,\theta_2)}\B = \ker \D \tilde P_{(\theta_1,\theta_2)}\oplus \ker \D G_{(\theta_1,\theta_2)}$ if and only if the determinant of the matrix $[\D G^T| \,\D\tilde P^T](\theta_1,\theta_2)$ is nonzero:
\begin{equation}
\T_{(\theta_1,\theta_2)}\B = \ker \D \tilde P_{(\theta_1,\theta_2)}\oplus \ker \D G_{(\theta_1,\theta_2)} \iff \det \begin{bmatrix}
\D G_{(\theta_1,\theta_2)}\\
\D \tilde P_{(\theta_1,\theta_2)}
\end{bmatrix}\not = 0.
\end{equation}
Examination of the matrix $[\D G^T| \,\D\tilde P^T](\theta_1,\theta_2)$ shows that $\det [\D G^T| \,\D\tilde P^T](\theta_1,\theta_2)$ is invariant under nonzero scaling of $(\theta_1,\theta_2)$.
I.e.,
\begin{equation}
\forall k \not = 0: 
\det \begin{bmatrix}
\D G_{(k\theta_1,k\theta_2)}\\
\D \tilde P_{(k\theta_1,k\theta_2)}
\end{bmatrix}
=
\det \begin{bmatrix}
\D G_{(\theta_1,\theta_2)}\\
\D \tilde P_{(\theta_1,\theta_2)}
\end{bmatrix}.
\end{equation}
In order to show that condition \eqref{eq:DG_conn} holds on $\B$, it therefore suffices to show that $\det [\D G^T| \,\D\tilde P^T](\theta_1,\theta_2)$ is nonzero whenever $\theta_1^2 + \theta_2^2 = 1$, or equivalently that
\begin{equation}
\forall 0 \leq \tau < 2\pi: \det \begin{bmatrix}
\D G_{(\cos\tau,\sin\tau)}\\
\D \tilde P_{(\cos\tau,\sin\tau)}
\end{bmatrix}\not = 0.
\end{equation}
This is indeed the case, as illustrated by the numerical proof offered in Figure \ref{fig:det_plot}.
Thus $\T \B = \ker \D G \oplus \ker \D P$ and therefore condition \eqref{eq:DG_conn} holds.

Denoting $\varphi:= \tilde P(\theta_1,\theta_2)$, the dynamics $g:\M \to \T \M$ are given by $g(P(\theta_1,\theta_2)) = g \circ F(\varphi) =  \D F(\varphi)\cdot 1 = (a_1\cos(\varphi-\delta),a_2\cos(\varphi+\delta))$,
or 
\begin{equation*}
g(P(\theta_1,\theta_2)) = s_h \frac{\left[\left(\frac{\theta_2}{a_2}-\frac{\theta_1}{a_1}\right)\frac{a_1\cos\delta}{2\sin\delta} + \left(\frac{\theta_2}{a_2}+\frac{\theta_1}{a_1}\right)\frac{a_1\sin\delta}{2\cos\delta}, \left(\frac{\theta_2}{a_2}-\frac{\theta_1}{a_1}\right)\frac{a_2\cos\delta}{2\sin\delta} - \left(\frac{\theta_2}{a_2}+\frac{\theta_1}{a_1}\right)\frac{a_2\sin\delta}{2\cos\delta} \right]^T}{D}.
\end{equation*} 

\begin{figure}[h]
	\centering 
	\includegraphics[width=1.0\textwidth]{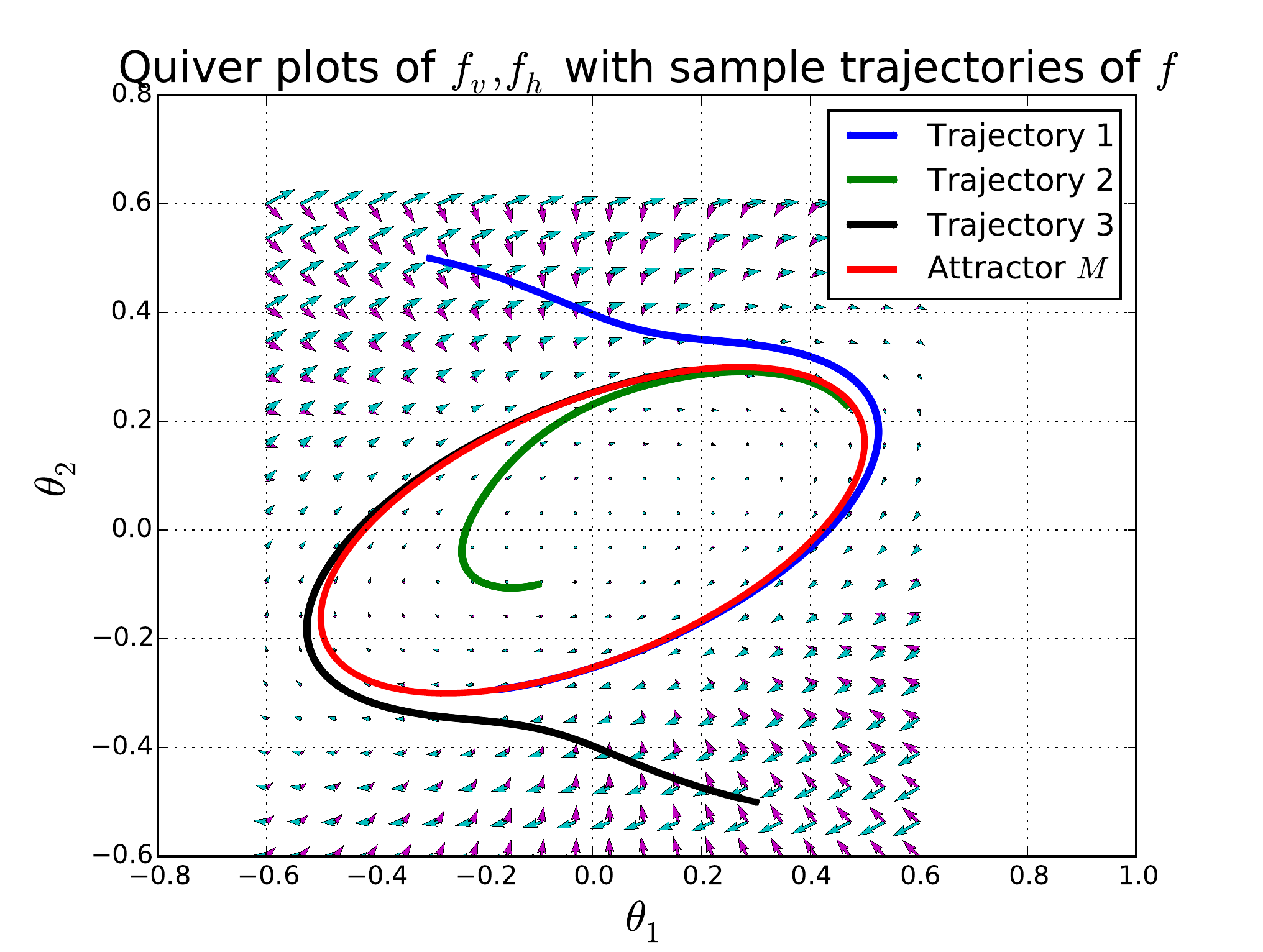}
	\caption[]{Quiver plots of $f_h$ (cyan) and $f_v$ (magenta) are shown together with three sample trajectories of $f:= f_h + f_v$.
	The attractor, $M$, is shown in red.} %
	\label{fig:all}
\end{figure}

We now have all of the ingredients necessary to compute 
$$f_h(\theta_1,\theta_2):= T_{(\theta_1,\theta_2)}^{-1} \left[T_{P(\theta_1,\theta_2)} \D P_{(\theta_1,\theta_2)}  T_{(\theta_1,\theta_2)}^{-1}\right]^\dagger T_{P(\theta_1,\theta_2)}g(P(\theta_1,\theta_2))$$
as in equation \eqref{eq:fh_def}, and 
$$f_v(\theta_1,\theta_2) = -\alpha(\theta_1,\theta_2)\Pi^P_{(\theta_1,\theta_2)} \nabla V(\theta_1,\theta_2)$$
as in equation \eqref{eq:fv_def}.
For the purpose of making Figure \ref{fig:all} easy to interpret, we chose 
\begin{align*}
a_1 &= 0.5 \\
a_2 &= 0.3\\
\delta &= 0.5\\
s_v &:= 5\times 10^{-3}\\
R &:= 10^{-4}\\
s_h&:= 0.75s_v\\
\alpha(x) &:= \frac{s_v R}{s_v + \|R \Pi^P_x\nabla V(x)\|}.
\end{align*}
Three sample trajectories of $f$ are shown in Figure \ref{fig:all}, along with the attractor $M$ and vector fields $f_h$ and $f_v$.

\subsubsection{The dynamic double pendulum}
In actual physical systems, one can only directly influence accelerations and not velocities through the application of force.
The vector field $f$ constructed in \S \ref{sec:kin_pendulum} is a direct application of the general construction in this paper, but assumes the ability to directly influence velocities.
If the construction in \S \ref{sec:construction} is applied directly to the phase space of a physical system, it will in general produce a non-physical vector field (one for which the derivative of position is not velocity).
One may therefore justifiably worry that the construction of \S \ref{sec:construction} is not applicable to physical systems.

However, there exist a variety of techniques for approximating the dynamics of ``first-order'' vector fields (in which velocities are directly influenced) by ``second-order'' vector fields \citep{revzen2012dynamical,koditschek1987adaptive}.
In order to apply the technique of \S \ref{sec:construction} to obtain physical vector fields on the phase spaces of physical systems, we thus propose the following general (intentionally vague) strategy.
First, define the quantities $\M \subseteq \B$, $P:\B\to \B$, $g:\M\to \T \M$ and construct the vector field $f:\B \to \T \B$ as in \S \ref{sec:construction}.
Second, define a ``second-order'' vector field whose dynamics approximate those of $f$.

\begin{figure}[h]
	\centering 
	\includegraphics[width=1.0\textwidth]{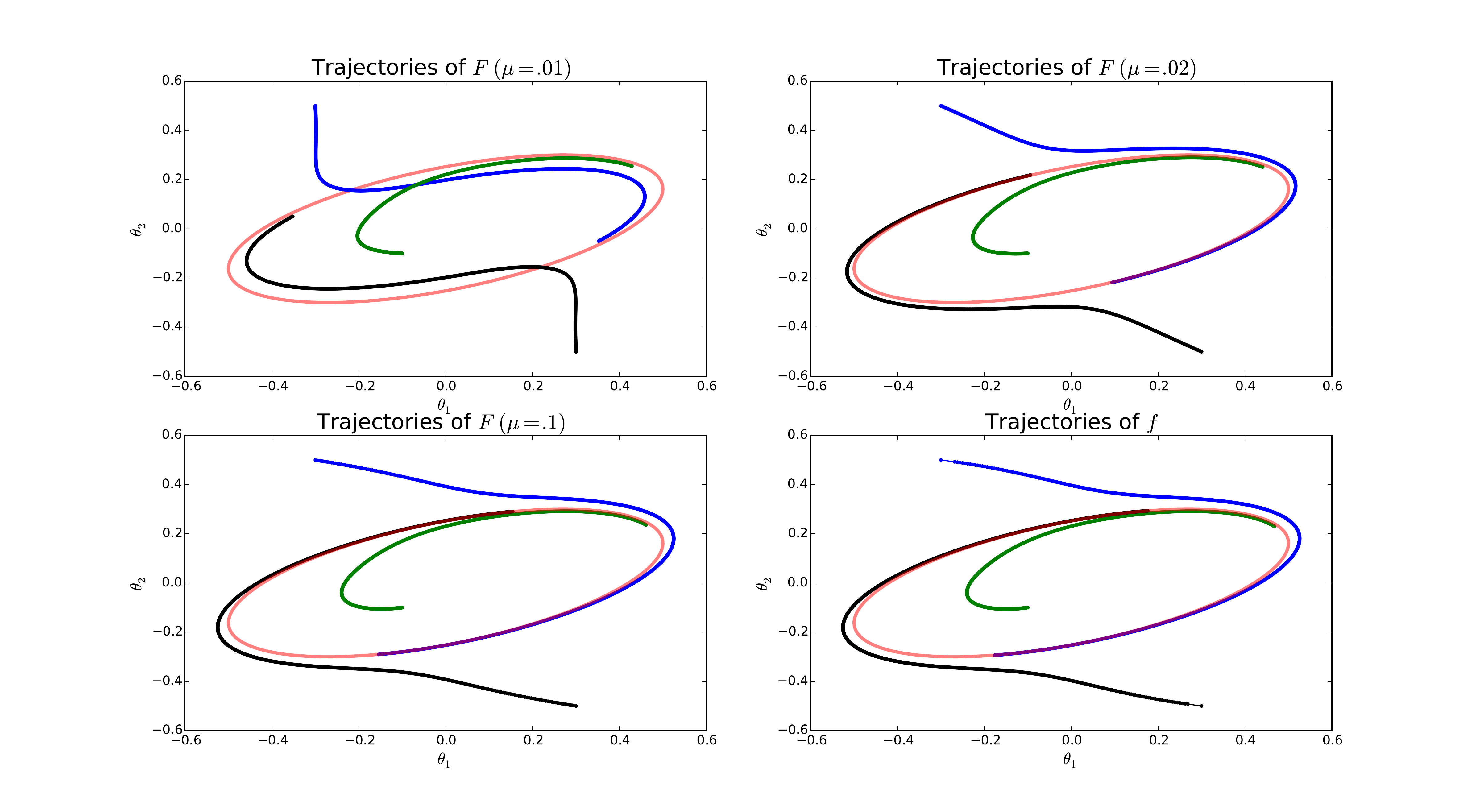}
	\caption[]{For various values of $\mu$, sample trajectories of $F$ having the same initial values (initial values of $\omega_1,\omega_2$ are all zero here) projected onto the $(\theta_1,\theta_2)$ plane are shown in the first three plots.
		Trajectories of $f$, the vector field defined in the example in \S \ref{sec:kin_pendulum}, are shown in the bottom-right plot.
		The attractor, $M$, is shown in red.
		These plots indicate that the dynamics of $f$ can be approximated arbitrarily well by increasing $\mu$.} %
	\label{fig:F_mu_subplots}
\end{figure}

We continue the example in \S \ref{sec:kin_pendulum} with a specific version of this strategy using a technique inspired by Koditschek (\citet{koditschek1987adaptive} \S 3.3).
For notational purposes, define $\theta:= (\theta_1,\theta_2)$ and $\omega:= (\omega_1,\omega_2)$.
For $\mu > 0$, define the vector field $F:\T \B \to \T(\T \B)$ by
\begin{equation}
F(\theta,\omega) = \left[\omega, \D f_\theta f(\theta) - \mu(\omega - f(\theta))\right]^T,
\end{equation}
 so that dynamics on $\T \B \cong \B \times \R^2$ are given by
\begin{align*}
\dot \theta &= \omega\\
\dot \omega &= \D f_\theta f(\theta) - \mu\left(\omega-f(\theta)\right).
\end{align*}
Ideas from the theory of normal hyperbolicity can likely be used to show that for sufficiently large $\mu$, each trajectory of $F$ asymptotically approaches approaches a corresponding trajectory in the invariant manifold $\{(\theta,\omega)|\omega = f(\theta)\}$.
However, noncompactness $\{(\theta,\omega)|\omega = f(\theta)\}$ complicates matters.
Such an analysis is outside the scope of this work, so we simply present numerical results in Figure \ref{fig:F_mu_subplots}.

\section{Lyapunov functions on $\M$ extend to $\B$}
\label{sec:lyap_extend}
Let $\eta:\M \to \R$ be any function.
Define the \concept{pullback} function $P^*\eta:\B\to \R$ by $P^*\eta:= \eta \circ P$.
Then since $P$ commutes with the flow, we have:
\begin{equation}\label{eq:lyap_pullback}
P^* \eta\circ \phi_t = \eta \circ P \circ \phi_t = \eta \circ \phi_t \circ P.
\end{equation}

\begin{Prop}
	Let $\eta:\M \to \R$ be a Lyapunov function for the dynamics $g$ on $\M$ corresponding to some invariant compact subset $K$ of $\M$. 
	Then the function $P^*\eta + V: \B \to \R$ is a Lyapunov function for $K$ on $\B$, where $P^*\eta:= \eta \circ P$ is the pullback of $\eta$.
\end{Prop}

\begin{proof}
	It follows from equation \eqref{eq:lyap_pullback} that $P^*\eta(\phi_t(x_0))$ is strictly decreasing if $x_0 \in \B \ct P^{-1}(K)$ and is identically zero if $x_0 \in P^{-1}(K)$.
	By assumption, $V(\phi_t(x_0))$ is strictly decreasing for all $x_0 \in \B \ct \M$ and is identically zero if $x_0 \in \M$. 
	It follows that $(P^*\eta + V)(\phi_t(x_0))$ is monotone decreasing for all $x_0 \in \B \ct K$ 
	
	Since $V$ is strictly positive on $\B \ct \M$ and zero on $\M$, and $P^*\eta$ is strictly positive on $\B \ct P^{-1}(K)$ and zero on $P^{-1}(K)$,
	it follows that $P^*\eta + V$ is strictly positive on $\B \ct K$ and is zero on $K$.
	Since $K$ is compact, the Lyapunov theorem (\citet{wilson1967structure} Theorem 3.1) implies that $K$ is asymptotically stable (though its basin of attraction might not be all of $\B$) and $P^*\eta + V$ is a Lyapunov function.
\end{proof}

\section{Further generalizing our methods}
\label{sec:generalizing}
\subsection{Topological limitations of the construction in \S \ref{sec:construction}}
The construction in \S \ref{sec:construction} produces dynamics on $\B$ such that $\M$ is an asymptotically stable invariant manifold with basin of attraction $\B$.
However, we assumed that $\M$ was a level set of a submersion $G:\B\to \R^{n-k}$.
Topological arguments show that this is not true of all embedded submanifolds, so our construction in \S \ref{sec:construction} is not applicable to all embedded submanifolds $\M$ of $\R^n$.
One way to fix this is to work completely within the abstract framework of fibered manifolds and fiber bundles and replace the function $G$ with a different connection playing the role of $\ker \D G$.
We plan to elucidate these ideas in a future publication.
\subsection{Non-compact attractors}
Our construction in \S \ref{sec:construction} can also be applied in the case that $\M$ is not compact.
However, there are additional technical conditions which must be considered in order to ensure completeness of the resulting vector field $f$, as well as asymptotic stability of $\M$.
Indeed, merely defining a notion of asymptotic stability of noncompact attractors is subtle -- it depends, in general, on the choice of distance function on $\R^n$ and can not be a purely topological notion -- see, e.g., the discussion following Theorem 3.3 of \citet{wilson1969smoothing}.
Additionally, both stating and proving a theorem like our Theorem \ref{th:persistence_of_our} for noncompact $\M$ becomes more involved.
Results on persistence of noncompact NHIMs have recently been proved by \citet{eldering2013normally}, but we know of no explicit results in the literature on existence and persistence of asymptotic phase for noncompact $\M$ (however, see the final paragraph on p. 4 of \citet{eldering2013normally}). 
All of these technical details take us rather far afield from the core ideas of this work; for this reason, we chose not to pursue them further and restricted ourselves to the case of constructing compact invariant manifolds.

\section{Discussion}
As a notion of controller design, ``anchoring a template'' seems to be a particularly powerful one.
It expresses the idea that a controller takes a complex anchor system and reduces its behavior to that of a simpler, better understood template system.
Normal hyperbolicity and the notion of unique asymptotic phase provide one natural way to express the template-anchor relationship.

What we have shown is that a broad class of NHIM based template-anchor systems can be reverse engineered -- i.e. they can be broken into mathematical building blocks, each of which contributes a clearly defined functionality, and put back together from those blocks.
Furthermore, given such blocks an embedded template manifold can be made normally hyperbolic and endowed with a nearly arbitrary choice of unique asymptotic phase, producing a system which is robust to perturbations and modeling errors.

The key insights enabling our construction are the following.
First, the explicit formula for constructing $f_h$ using standard matrix computations was the creation of a coordinate change $T_x$, as defined in Equation \eqref{eq:T_def}, under which the horizontal and vertical components of the flow become orthogonal.
Orthogonality in this new coordinate system enables the lift $f_h$ to be constructed via the standard Moore-Penrose pseudoinverse while ensuring that $\M$ is an invariant manifold of $f_h$. 
The second insight is that if $f_v$ takes values in the vertical bundle $\ker \D P$, then $f_v$ and $f_h$ don't interfere with each other.
Thus if $f_v$ stabilizes $\M$, $f_h$ and $f_v$ can be combined to give a vector field $f$ rendering $\M$ asymptotically stable with asymptotic phase, anchoring template dynamics, and enabling a broad range of applications.

\subsection{Acknowledgements}
This work was supported by ARO Morphologically Modulated Dynamics \#W911NF-14-1-0573 to S. Revzen.
The authors wish to thank Jessy Grizzle, Ralf Spatzier, Jaap Eldering, and George Council for helpful conversations.

\appendix
\section{Moore-Penrose pseudoinverse}
\label{app:pseudoinverse}
Let $V$ and $W$ be real inner product spaces, and let $L:V \to W$ be a linear map.
Let $L^T: W \to V$ denote the \concept{transpose} or \concept{adjoint} of the linear map $L$; that is, $L^T$ is the unique linear map such that for all $v \in V$ and $w \in W$, $\langle Lv,w \rangle = \langle v,L^Tw \rangle $.
The transpose of a composition of linear maps satisfies $(AB)^T = B^TA^T$ from which it follows that, e.g., $(ABC)^T = C^TB^TA^T$.
If $V$ and $W$ are Euclidean spaces with the Euclidean inner product and $L$ is identified with its matrix representation with respect to the standard Euclidean bases, then $L^T$ is the ordinary matrix transpose.

A \concept{Moore-Penrose pseudoinverse}, or more succinctly \concept{pseudoinverse}, of $L$ is a linear map $L^\dagger:W \to V$ satisfying the following four properties

\begin{align*}
1. \,\, L L^\dagger L = L \qquad\ 2.\,\,L^\dagger L  L^\dagger = L^\dagger \qquad 3.\,\, (L L^\dagger)^T=LL^\dagger \qquad 4.\,\, (L^\dagger L)^T = L^\dagger L.
\end{align*}

Given any linear map $L:V\to W$, a linear map $L^\dagger$ satisfying the above properties exists and is unique \citep{penrose1955generalized}. 
Note that the definition of the pseudoinverse depends entirely on the choice of inner products for $V$ and $W$ because the definition of the adjoint depends on the inner products; a different choice of inner products would result in a different pseudoinverse.

It is useful to think about the pseudoinverse in terms of orthogonal projections.
It can be shown that $P = LL^\dagger: W\to W$ is orthogonal projection onto $L(V)$, and $Q = L^\dagger L: V \to V$ is orthogonal projection onto $(\ker L)^\perp$.

We will be concerned with the case in which $L$ is surjective. 
In this case, it can easily be shown that $LL^T$ is invertible.
Combining the first and fourth properties above shows that $L = L(L^\dagger L)^T = L L^T (L^\dagger)^T$.
Taking the adjoint of this equation shows that
$L^T = L^\dagger L L^T$, so 
\begin{equation}\label{eq:pseudoinverse}
L^\dagger = L^T (LL^T)^{-1}
\end{equation}
since $LL^T$ is invertible in the case that $L$ is surjective. 
It follows that, in the case that $L$ is surjective, $LL^\dagger:W\to W$ is the identity map.

\section{Differential topology} \label{app:differential_topology}
% the \\ insures the section title is centered below the phrase: AppendixA

If $F:X \to Y$ is a map between sets and $U\subseteq X$, then $F|_{U}:U \to Y$ denotes the \concept{restriction} of $F$ to $U$.
Given any subset $V \subseteq Y$, $F^{-1}(V):= \{x \in X| F(x)\in V\}$ is the \concept{pre-image} of $V$ under $F$.

If $U$ is any subset of a topological space $X$, we let $\mathring U$ denote its interior, $\bar U$ denote its closure, and $\partial U$ denote its boundary.

If $M$ is an $m$-dimensional $\C^r$ manifold ($r \geq 1$) and $p \in M$, let $\T_p M$ denote the tangent space to $M$ at $p$.
We recall that the tangent bundle $\T M$ of $M$ has a natural topology and smooth structure making $\T M$ into a $\C^{r-1}$ manifold of dimension $2m$.
If $F:M \to N$ is a $\C^r$ map between $\C^r$ manifolds and $x \in M$, we denote by $\D F_x:\T_x M\to \T_{F(x)}N$ the \concept{differential} of $F$ at $x$, which is a linear map.
We recall that the map $\D F: \T M \to \T N$ defined by $\D F(x,v) = \D F_x v$ is a $\C^{r-1}$ map.
$F$ is an immersion if $\D F_x$ is injective at each $x \in M$.
$F$ is a submersion if $\D F_x$ is surjective at each $x \in M$.
$F$ is a $\C^r$ local diffeomorphism if for each $x\in M$, there exists a neighborhood $U$ containing $x$ such that $F:U \to F(U)$ is a diffeomorphism.
$F$ is a $\C^r$ local diffeomorphism if and only if $F$ is $\C^r$ and $F$ is both an immersion and submersion.

A $\C^r$ submersion $F:M \to N$ between $m$ and $n$-dimensional manifolds has the special property that for any $x \in M$, there exist open sets $U \ni x$ and $V \ni F(x)$ together with $\C^r$ diffeomorphisms $\varphi:U\to \hat U \subset \R^m$ and $\psi:V \to \hat V\subset \R^n$ such that 
\begin{equation}
\psi \circ F \circ \varphi^{-1}(x_1,\ldots,x_n,\ldots,x_m) = (x_1,\ldots,x_n,0,\ldots,0).
\end{equation}

A map $F:U \to Y$ between topological spaces is a \concept{topological embedding} if it is a homeomorphism onto its image. 
If $F:U \to Y$ is a $\C^r$ map ($r \geq 1$), 
A $\C^r$ map $F:U\to Y$ between $\C^r$ manifolds is a \concept{$\C^r$ embedding} if it is a topological embedding and an immersion. 

If $M$ is any $\C^r$ $m$-dimensional manifold ($r \geq 0$), a subspace $K$ of $\B$ with the subspace topology is a $\C^r$ $k$-dimensional \concept{embedded submanifold} if it is a $\C^r$ manifold and the inclusion map $i:K \hookrightarrow M$ is a $\C^r$ embedding.
A \concept{proper map} between topological spaces is a map for which the pre-image of any compact set is a compact set.
A submanifold is \concept{properly embedded} if it is an embedded submanifold and additionally the inclusion map is a proper map; equivalently, a submanifold is properly embedded if and only if it is an embedded submanifold which is also a closed subset of the ambient manifold.

\section{Fibered manifolds, fiber bundles, and connections}
\label{app:fiber_connections}
A triple $(\B,\tilde P,\tilde \M)$, where $\tilde P$ is a $\C^r (r\geq 1)$ surjective submersion $\tilde P:\B \to \tilde \M$ between $\C^r$ manifolds $\B$ and $\tilde \M$, is called a \concept{fibered manifold} (\citet{kolar1999natural} \S 2).
$\B$ is called the \concept{total space}, and $\tilde \M$ is called the \concept{base}.
We sometimes simply refer to $\tilde P:\B \to \tilde \M$ as being a fibered manifold or even more succinctly we may just refer to $\B$ as being a fibered manifold.
Because $\tilde P$ is a submersion, a fibered manifold has the property that for any $x \in \B$, there exist $\C^r$ \concept{fiber charts} $(U,\varphi)$ and $(V,\psi)$ with $U \ni x$ and $V \ni \tilde P(x)$ open sets and $\varphi:U\to\hat U \subset \R^n$ and $\psi:V\to\hat V \subset \R^m$ diffeomorphisms such that:

\begin{equation}
\psi \circ \tilde P \circ \varphi^{-1}(x_1,\ldots,x_n,\ldots,x_m) = (x_1,\ldots,x_n,0,\ldots,0).
\end{equation}

We define a $\C^0$ fibered manifold to be a triple $(\B,\tilde P,\tilde \M)$ to be a continuous $\tilde P:\B \to \tilde \M$ to be a continuous map between topological manifolds so that given any $x \in \tilde \M$ there exist $\C^0$ fiber charts containing $x$ and $\tilde P(x)$ as above.

A $\C^r (r \geq 0)$ \concept{fiber bundle} (\citet{kolar1999natural} \S 9) is a tuple $(\B,\tilde P,\tilde \M,F)$ with $F$ a $\C^r$ manifold such that $(\B,\tilde P,\tilde \M)$ is a $\C^r$ fibered manifold and additionally for each $x \in \tilde \M$ there exists an open set $U \subseteq \tilde \M$ containing $x$ such that $\tilde P^{-1}(U)$ is $\C^r$ diffeomorphic to $U \times F$ via a diffeomorphism which respects fibers:

\begin{equation}\label{eq:fiber_bundle_diagram}
\begin{tikzcd}
&\tilde P^{-1}(U) \arrow{rr}{\rho}\arrow{dr}{\tilde P} & &U \times F\arrow{dl}{\text{pr}_1}\\
& &U
\end{tikzcd}
\end{equation}
where $\text{pr}_1:U \times F \to U$ is projection onto the first factor.
The map $\rho$ is referred to as a \concept{local trivialization}.
A fiber bundle $\B,\tilde P,\tilde \M,F)$ is \concept{trivial} if there exists a local trivializtion $\rho: \B \to \tilde \M \times F$.
A $\C^r$ map $\sigma:\tilde \M \to \B$ such that $\tilde P \circ \sigma$ is the identity map on $\tilde \M$ is called a $\C^r$ \concept{section} of $\tilde P$.
Given an open set $U \subseteq \tilde \M$, a $\C^r$ map $\sigma:U \to \B$ such that $\tilde P \circ \sigma$ is the identity map on $U$ is called a $\C^r$ \concept{local section} of $\tilde P$.

A (real) $\C^r$ \concept{vector bundle} of rank $k$ is a $\C^r$ fiber bundle such that each fiber of $\tilde P$ is endowed with the structure of a real $k$-dimensional vector space and such that any $x$ in the base space has a neighborhood $U$ and local trivialization $\rho: \tilde P^{-1}(U)\to U \times F$ such that the restriction of $\rho$ to each fiber of $\tilde P$ is a linear vector space isomorphism (\citet{lee2013smooth} Chapter 10).
Examples of $\C^{r-1}$ vector bundles include the tangent bundle and normal bundle of a $\C^r$ manifold, where $r \geq 1$.
The \concept{zero section} of a vector bundle is a section $\sigma$ sending each point $p \in \M$ to the zero vector in $P^{-1}(p)$.
We will also use the term $\concept{zero section}$ to refer the image $\sigma(\tilde \M)$ of the zero section $\sigma$.

Given $0 \leq j \leq r$, $1\leq m \leq k$, a $\C^j$ rank $m$ \concept{subbundle} $(\B_s,\tilde P_s,\tilde \M,F_s)$ of a $\C^r$ rank $k$ vector bundle $(\B,\tilde P,\tilde \M,F)$ is a $\C^j$ rank $m$ vector bundle in which $\B_s$ is a $\C^j$ embedded submanifold of $\B$, $\tilde P_s = \tilde P|_{\B_s}$, each fiber $\tilde P_s^{-1}(x) = \tilde P^{-1}(x)\cap \B_s$ is a linear subspace of $\tilde P^{-1}(x)$, and the vector space structure on $\tilde P_s^{-1}(x)$ is the vector space structure inherited as a subspace of $\tilde P^{-1}(x)$.
In practice, the following \concept{local frame criterion} is often easier to check: $(\B_s,\tilde P_s,\tilde \M,F_s)$ is a $\C^j$ rank $m$ subbundle of $(\B,\tilde P,\tilde \M,F)$ if and only if for every point $x \in \tilde \M$, there is a neighborhood $U\subseteq \tilde \M$ containing $x$ and $\C^j$ local sections $\sigma_1,\ldots,\sigma_m:U\to \B$ of $\tilde P:\B \to \tilde \M$ such that for all $y \in U$: $\sigma_1(y),\ldots, \sigma_m(y)$ form a basis for the vector space $\tilde P_s^{-1}(y)$ (\citet{lee2013smooth} Lemma 10.32). 

Let $ r \geq 1$.
Given a $\C^r$ fibered manifold $(\B,\tilde P,\tilde \M)$, for any $x \in \B$ we define the \concept{vertical space} $\Ver_x\B:= \ker \D P_x$. 
The \concept{vertical bundle} $\Ver \B \subseteq \T \B$  is the union of all of the vertical spaces (i.e., $\Ver \B := \bigcup_{x\in\B}\Ver_x\B$),
endowed with the unique topology and smooth structure 
making the vertical bundle into a $\C^{r-1}$ subbundle  of the $\C^{r-1}$ tangent bundle.
A $\C^{r-1}$ \concept{connection} for $\tilde P$   is a $\C^{r-1}$ subbundle $\Hor \B$ of $\T\B$ such that for each $x \in \B$, $\T_x \B = \Hor_x\B \oplus \Ver_x\B$.
We also refer to $\Hor \B$ as a \concept{horizontal bundle}. 
To every connection $\Hor \B$ corresponds a vertical-valued projection $\pi: \T \B \to \Ver \B$ such that for any $x \in \M$, $\pi|_{\T_x \B}$ is linear projection onto $\Ver_x \B$ with kernel $\Hor_x\B$.

Note that since every fiber bundle (with manifold base and total space) is also a fibered manifold, the definitions of the preceding paragraph apply to fiber bundles.
A $\C^{r-1}$ connection is \concept{complete} if, given any $\C^{r-1}$ path $\gamma:[0,1]\to \M$ and $x \in P^{-1}(\gamma(0))$, there exists a $\C^{r-1}$ \concept{lift} $\tilde \gamma:[0,1]\to \B$ satisfying $P \circ \tilde \gamma = \gamma$, $\tilde \gamma(0) = x$, and $\forall t \in [0,1]: \dot {\tilde \gamma}(t) \in \Hor_{\tilde \gamma(t)}\B$.
For $r \geq 2$, using results of \citep{del2015complete,kolar1999natural} together with results of approximation theory (\citet{lee2013smooth} Chapter 6), it can be shown that a $\C^r$ fibered manifold is a fiber bundle if and only if it admits a complete $\C^{r-1}$ connection.

Let $s \geq r \geq 0$. 
We will use the phrases \concept{$\C^r$ fibered manifold with $\C^s$ fiber} and \concept{$\C^r$ fiber bundle with $\C^s$ fiber} to refer to $\C^r$ fibered manifolds and fiber bundles whose individual fibers are $\C^s$ embedded submanifolds.  

\section{Dynamical systems theory}
\label{app:dynamical_systems}
Given a $\C^r$ ($r \geq 1$) manifold $M$, a $\C^{r-1}$ map $F:M \to \T M$ is called a \concept{$\C^{r-1}$ vector field} if $\pi \circ F:M \to M$ is the identity map, where $\pi: \T M \to M$ is the natural projection.
To use the vocabulary introduced in Appendix \ref{app:fiber_connections}, a $\C^{r-1}$ vector field $F:\M \to \T \M$ is a $\C^{r-1}$ section of the vector bundle projection $\pi:\T\M \to \M$.
A \concept{trajectory}, \concept{solution}, or \concept{integral curve} of $F$ is a $\C^r$ curve $\gamma : J\to M$ such that $\dot \gamma(t) = F(\gamma(t))$ for every $t \in J$, where $J\subseteq \R$ is an interval.
A trajectory is \concept{maximal} if its domain $J$ cannot be extended to any larger interval.
Under appropriate conditions, maximal integral curves exist, are unique, and the vector field $F$ admits an associated \concept{flow} $\phi:W \to M$, where the \concept{maximal flow domain} $W\subseteq \R \times M$ is an open set (see, e.g., Chapter 8 of \citet{odeHirschSmale} or Chapter 9 of \citet{lee2013smooth}).
We sometimes write $\phi_t(x):= \phi(t,x)$.
For all $x \in M$ and $t,s \in \R$ for which the following expression is defined, the flow satisfies the \concept{group properties} $\phi_t\circ \phi_s(x) = \phi_{t+s}(x)$ and $\phi_0(x)=x$.
For all $(t,x) \in W$, the flow also satisfies $\frac{\partial}{\partial t}\phi_t(x) = F(\phi_t(x))$,
so that for any fixed $x \in M$ the curve $t \mapsto \phi_t(x)$ is an integral curve through $x$.
The maximal flow domain $W$ is defined so that the restriction of $\phi$ to any set of the form $W \cap \{p\} \times M$ is a maximal integral curve.

\section{Connection lemmas}
\label{app:connection_lemmas}
We will use the following lemma on existence and regularity of connections on fibered manifolds in which the base space is an embedded $k$-dimensional submanifold of the total space, which in turn is an open subset of $\R^n$.
\begin{Lem}\label{lem:connection}
	Let $\B$ be an open subset of $\R^n$ and let $\M\subseteq \B$ be a compact $\C^r$ $k$-dimensional embedded submanifold.
	Let $P:\B \to \M$ be a $\C^1$ map such that $(\B,P,\M)$ is a $\C^1$ fibered manifold.
	For $x \in \B$, let $\Ver \B = \ker \D P$ be the $\C^0$ vertical bundle $\Ver \B$.
	Then there exists a $\C^{r-1}$ connection (so that $\Hor \B$ $\T \R^n = \Ver \B \oplus \Hor \B)$ such that $\Hor\B|_\M = \T \M$.
\end{Lem} 

\begin{proof}
	We first show that we can define such a $\C^0$ connection $\Hor U$ on some neighborhood $U$ of $\M$ with the property that $\forall x \in \M: \Hor_x\B = \T_x\M$. 
	
	$\M$ is a $\C^r$ embedded submanifold of $\B$ and $\T \M$ is a $\C^{r-1}$ subbundle of $\T \B$.
	This implies that for any $x \in \M$, there exists an open neighborhood $U_x \ni x$ and a $\C^{r-1}$ diffeomorphism $\psi:\T\B|_{U_x} \to U_1 \times U_2 \subseteq \R^n \times \R^n$ such that $\psi(\T\M|_{U_x}) = \{(y_1,\ldots,y_k,\ldots,y_n,v_1,\ldots,v_k,\ldots,v_n)\in \R^k \times \R^n|y_{k+1}=\ldots=y_{n} = v_{k+1}=\ldots=v_n=0\}$.
	Define a $\C^{r-1}$ connection $\Hor U_x$ by 
	$$\Hor_y U_x := (\D\psi_y)^{-1}\left(\{\psi(y)\}\times \{(v_1,\ldots,v_k,\ldots,v_n)\in \R^n| v_{k+1}=\ldots=v_n=0\}\right).$$
	Let $\pi^{U_x}:\T\B_{U_x}\to\Ver\B$ denote the vertical-valued projection with kernel equal to $\Hor U_x$, and note that $\pi^{U_x}$ is a $\C^0$ map since $\Ver \B$ is a $\C^0$ subbundle of $\T\B$\footnote{This follows by, e.g., repeating the proof of (\citet{lee2013smooth} Theorem 10.34), replacing ``smooth'' everywhere with ``continuous.''}.
	Let $\{\varphi^x\}_{x\in\M}$ be a partition of unity subordinate to the open cover $\{U_x\}_{x\in\M}$, let $U:= \bigcup_{x\in\M}U_x$, and define the $\C^0$ vertical-valued map $\pi^U$ by $$\pi^U:= \sum_{x\in\M}\varphi^x\pi^{U_x}.$$
	$\pi^U$ is also a vertical-valued projection since for each $x\in\M$ this sum is a convex combination of the vertical-valued projections $\pi^{U_x}$, and projections with common image are closed under convex combinations.
	Since $\pi^U$ is $\C^0$, the connection $\Hor U$ it defines is also $\C^0$.
	Note that for any $x \in \M$ with $x \in U_y \cap U_z$, $\ker \pi^{U_y}_x = \ker \pi^{U_z}_x = \T_x\M$, and hence $\pi^{U_y}_x = \pi^{U_z}_x$.
	It follows that $\ker \pi^U|_\M = \T \M$.

	We next define a $\C^0$ connection on the open set $W:= \B \ct \M$.
	We simply define the horizontal space $\Hor W:= (\Ver\B|_W)^\perp$.
	Let $\pi^W$ be the corresponding $\C^0$ vertical-valued projection.
	
	Next, let $\varphi_W, \varphi_U$ be a partition of unity subordinate to $W,U$. 
	We define $\pi':= \varphi_W \pi_W + \varphi_U \pi_U$.
	Once again using the fact that convex combination of any two linear projection operators with common image is again a projection onto the same image, it follows that $\pi'$ is a $\C^0$ vertical-valued projection on $\B$.
	Note that, since the support of $\varphi_W$ is contained in the complement of $\B$, $\ker \pi'|_\M = \T \M$.
	Letting $\Hor'\B$ denote the $\C^0$ connection corresponding to $\pi'$, it follows that $\Hor'|_\M$ is $\C^{r-1}$.
	Denote by $\Hor'\B$ the $\C^0$ horizontal bundle corresponding to $\pi'$.
	
	% NOTE: This is the only part of the proof where I used the assumption that everything is a subset of \R^n. Everything would generalize if I could figure out a different way to approximate which doesn't use identification of tangent spaces with \R^n.
	
	To complete the proof, we will use approximation techniques to approximate the $\C^0$ connection $\Hor'\B$ by a $\C^{r-1}$ connection. 
	Let $G_k$ be the $\C^\infty$ Grassmann manifold of $k$-dimensional linear subspaces of $\R^n$.
	By the Whitney Embedding Theorem (\citet{lee2013smooth} Chapter 6), we may consider $G_k$ to be an embedded submanifold of $\R^N$ for some $N > 0$, so given $S_1,S_2 \in G_k$ we may define $\|S_1-S_2\|$ to be the distance between $S_1$ and $S_2$ using the Euclidean norm on $\R^N$.
	We define a $\C^0$ map $J':\B \to G_k$ by $J'(x) = \Hor'_x\B$, where $\Hor'_x\B$ is viewed as a linear subspace of $\R^n$ after the standard natural identification of $\T_x\R^n$ with $\R^n$.
	Note that $J'|_\M$ is $\C^{r-1}$.
	
	We show that we can define a $\C^{r-1}$ map $J: \B \to G_k$ such for each $x \in \B$, $\T\R^n = J(x)\oplus\Ver_x\B$ and such that for each $x \in \M$, $J(x) = \T_x\M$.
	For each $x \in \B$, define  $\delta(x):=\sup\{r > 0| \|J'(x)-S\|<r\implies \T_x\R^n = S \oplus J'(x)\}$.
	Since $\Ver\B$ is a $\C^0$ vector bundle, it can be shown that $\delta$ is continuous.
	The Whitney Approximation Theorem (\citet{lee2013smooth} Chapter 6) shows that there exists a $\C^{r-1}$ map $J:\B\to G_k$ such that $J|_\M = J'|_\M$ and for all $x \in \B$, $\|J(x)-J'(x)\| \leq \delta(x)$.
	It follows that $\T_x\R^n = J(x)\oplus \Ver_x\B$ for each $x \in \B$.
	
	Taking $\Hor_x \B = J(x)$ for each $x \in \B$ completes the proof.
\end{proof}

\begin{Lem}\label{lem:hor_lift_smooth}
	Let $(\B,P,\M)$ be a $\C^r$ fibered manifold. 
	The horizontal lift of a $\C^{r-1}$ vector field on $\M$ via a $\C^{r-1}$ connection is $\C^{r-1}$.
\end{Lem}

\begin{proof}
	Let $x \in \M$. Since $\Hor \B$ is a $\C^{r-1}$ subbundle of $\T\M$, there exist $k$ pointwise linearly independent $\C^{r-1}$ vector fields $v_1,\ldots,v_k,\ldots,v_n$ defined on a neighborhood $U$ of $x$ such that $v_1(y),\ldots,v_k(y)$ span $\Hor_y \B$ at each $y \in U$, and $v_{k+1}(y),\ldots,v_n(y)$ span $\Ver_y\B = \ker \D P_y$ at each $y \in U$.
	Since $\D P_y$ is full rank as a map into $\T_y \M$, 
	$\D P_yv_1(y),\ldots,\D P_yv_k(y)$ form a basis for $\T_y\M$. 
	We may complete this to a basis  
	\\$\D P_yv_1(y),\ldots,\D P_yv_k(y),w_{k+1}(P(y)),\ldots w_n(P(y))$ for $\T_y\B$.
	For $y \in U$, the matrix of $\D P_y$ with respect to the bases $v_1(y),\ldots,v_n(y)$ of $\T_y\B$ and $\D P_y v_1(y),\ldots, \D P_yv_k(y),w_{k+1}(P(y)),\ldots w_n(P(y))$ is of the form
	
	\begin{equation*}
	[\D P_y]=
	\left[
	\begin{array}{c|c}
	I & 0 \\
	\hline
	0 & 0
	\end{array}
	\right].
	\end{equation*}
	
	Since $P$ is a submersion, $P(U)$ is an open set.
    Any $\C^{r-1}$ vector field $g$ on $\M$ may be written on $P(U)$ as $g(y)=\sum_{i=1}^k \kappa_i(y)\D P_y v_i(y)$, for some $\C^{r-1}$ functions $\kappa_i:P(U)\subseteq \M \to \R$.
	It follows that on $U$, the lift $f_h$ of $g$ is given by $f_h(y) = \sum_{i=1}^k\kappa_i(P(y))v_i(y)$.
	$\kappa_i\circ P$ is a $\C^{r-1}$ function on $U$ for each $i$ wince $\kappa_i$ and $P$ are each at least $\C^{r-1}$, and all of the $v_i$ are $\C^{r-1}$ vector fields so it follows that $f_h$ is $\C^{r-1}$ on $U$.
	Since we have shown that $f_h$ is $\C^{r-1}$ on a neigborhood of each point, it follows that $f_h$ is $\C^{r-1}$.
\end{proof}

\section{Normally hyperbolic invariant manifolds}
\label{app:NHIMs}
In this section, we summarize some of the main results on normally hyperbolic invariant manifolds (NHIMs) which served as the motivation for our construction in \S \ref{sec:construction}.
For simplicity, we only consider NHIMs which are embedded submanifolds of Euclidean space.

Normally hyperbolic invariant manifolds (NHIMs) are generalizations of hyperbolic fixed points and periodic orbits.
Much of the theory of compact NHIMs was independently developed in the 1970s by Fenichel \citep{fenichel1971persistence,fenichel1973asymptotic,fenichel1977asymptotic} and Hirsch, Pugh, and Shub \citep{hirsch1977}.
Eldering has recently extended many of these results to the noncompact setting \citep{eldering2013normally}.
We only need results on compact NHIMs; we choose to follow Fenichel's treatment here.
Since we are interested only in the case of asymptotically stable invariant manifolds, we will define a special case of normal hyperbolicity which is suitable for our needs. 

Let $\M$ be a compact embedded $\C^r$ submanifold of $\R^n$, invariant under the flow $\phi_t(\cdot)$ defined on some neighborhood of $\M$.
Let $\Nor \M$ be the normal bundle of $\M$, and let $\Pi^N:\T\R^n|_\M \to \Nor \M$ be the family of linear projections such that $\Pi^N_x$ orthogonally projects each tangent space $\T_x \R^n$ onto $\Nor_x \M$, for each $x \in \M$.
We will suppress the subscript $x$ in much of the sequel when the notation becomes cumbersome unless we wish to emphasize the role of $x$.
We define the linear maps $A_t(p):\T_p\M \to \T_{\phi_{-t}(p)\M}$ and $B_t(p):\Nor_{\phi_{-t}(p)}\M \to \Nor_p\M$:
\begin{align*}
A_t(p) &:= \D \phi_{-t}(p)|_{\T_p \M}\\
B_t(p) &:= \Pi^{\Nor} \D \phi_{t}(\phi_{-t}(p))|_{\Nor_{\phi_{-t}(p)}\M}\\
\end{align*}
\begin{Def}
	(Generalized Lyapunov type numbers). We define the following \concept{generalized Lyapunov-type numbers} for each $p \in \M$:
	
	\begin{align*}
	\nu(p)&:= \inf\left\{a > 0| \lim_{t\to\infty} \|B_t(p)\|/a^t = 0\right\}\\
	\sigma(p)&:= \inf\left\{s \in \R|\lim_{t\to\infty}\|A_t(p)\|\|B_t(p)\|^s = 0 \right\}\\
	\tau(p)&:= \inf\left\{s\in \R| \lim_{t\to\infty}\|A_{-t}(p)\| \left(\| B_t(\phi_t(p))\| A_t(\phi_t(p))\|\right)^s = 0\right\}.
	\end{align*}
\end{Def}
While we use the Euclidean norm here, it can be shown that the values of the generalized Lyapunov-type numbers are independent of the choice of inner product\footnote{in fact, independent of the choice of any Riemannian metric on $\R^n$.} (inducing a norm) on $\R^n$.

\begin{Def}
	We say that $\M$ is \concept{stable} if, given any open neighborhood $U \supseteq \M$, there exists an open neighborhood $V \supseteq \M$ such that $\forall t >0: \phi_t(V) \subseteq V$.
	We say that $\M$ is \concept{asymptotically stable} if $\M$ is stable and there exists a neighborhood $W \supseteq \M$ such that $\forall p \in W, t > 0: \dist{(\phi_t(p))}{\M} = 0$.
	We say that $\M$ is \concept{exponentially stable} if $\M$ is asymptotically stable and furthermore there exists $C, T, \mu > 0$ such that (possibly after shrinking $W$) $\forall p \in W, t > T: \dist{(\phi_t(p))}{\M} \leq Ce^{-\mu t}$.
\end{Def}

\begin{Prop}\label{prop:stable_NHIMS_are_exp_stable}
	If $\nu(p) < 1$ for every $p\in \M$, then $\M$ is exponentially stable.
\end{Prop}

\begin{Def}
	(Asymptotically stable normally hyperbolic invariant manifolds). Let $\M$ be a compact invariant $\C^r$ submanifold of $\R^n$, invariant under the $\C^r$ flow $\phi_t(\cdot)$ defined on some neighborhood of $\M$.
	We say that $\M$ is \concept{$r$-normally hyperbolic} if for all $p \in \M: \nu(p)<1$ and $\sigma(p)<\frac{1}{r}$.
	Without further qualification, \concept{normally hyperbolic} will be taken to mean $1$-normally hyperbolic.
\end{Def}

\begin{Rem}
	It can be shown that this definition is equivalent to \concept{eventual relative $r$-normal hyperbolicity} in \citet{hirsch1977} in the case that the NHIM is asymptotically stable; this is because $\M$ is compact.
	This definition is weaker than the often-used \concept{immediate relative $r$-normal hyperbolicity} found in \citet{hirsch1977}, yet most of the same main results hold.
	Because \citet{eldering2013normally} considers noncompact NHIMs, the definition of normal hyperbolicity chosen in \citet{eldering2013normally} is equivalent to \concept{eventual absolute $r$-normal hyperbolicity} as defined in \citet{hirsch1977}, which is also stronger than our definition. 
\end{Rem}

For convenience, we restate Definition \ref{def:asymptotic_phase} from \S \ref{sec:intro} and Definition \ref{def:unique_asymptotic_phase} from \S \ref{sec:motivation} here.

\begin{Def*}
	We say that $\M$ has \concept{asymptotic phase} if for any $x \in \B$, there exists a unique $P(x)\in \M$ such that $$\lim_{t\to\infty}\|\phi_t(x)-\phi_t(P(x))\| = 0.$$
	We say that $\M$ has \concept{unique asymptotic phase} if $\M$ has the asymptotic phase and additionally for any $x \in \B$ and any $q \in \M$ not equal to $P(x)$, 
	$$\lim_{t\to\infty}\frac{\|\phi_t(x)-\phi_t(P(x))\|}{\|\phi_t(x)-\phi_t(q)\|}=0.$$
	We refer to $P$ as the \concept{phase map} or simply as \concept{phase}, and say that $\M$ has \concept{$\C^k$ unique asymptotic phase} if the map $P:\B\to \M$ is $\C^k$.
\end{Def*}

The following restatement of Proposition \ref{prop:phase} in \S \ref{sec:motivation} is a combination of results from \citet{fenichel1973asymptotic,fenichel1977asymptotic}, and Theorem 4.1 of \citet{hirsch1977}.

\begin{Prop-hand}{1} %\label{prop:phase}
	Let $\M \subseteq \R^n$ be a compact $\C^r$ $k$-dimensional embedded submanifold of $\R^n$, invariant under the flow $\phi_t(\cdot)$ of the vector field $f:\Q\to \T \Q$ defined on an open neighborhood $\Q \subset \R^n$ of $\M$.
	Assume that for all $p \in \M$, $\sigma(p)<1$ and $\nu(p)<1$. 
	Then the following holds:

    \begin{enumerate}
    	\item The stability basin $\B$ of $\M$ is invariantly fibered by $\C^r$ manifolds $W_q$.
    	Explicitly, $\phi_t(W_q) = W_{\phi_t(q)}$, and the collection $\{W_q\}_{q\in\M}$ is a partition of $\B$.
    	Each $W_q$ is $\C^r$ diffeomorphic to $\R^{n-k}$.
    	Each $W_q$ intersects $\M$ transversally in the point $q$.
    	
		\item Let $P:\B \to \M$ be the map that sends $x \in \B$ to $q$, where $x \in W_q$. 
		Then $P$ is a continuous map, and $(\B,P,\M)$ is a $\C^0$ fibered manifold with $(n-k)$-dimensional Euclidean fibers.
		
		\item Let $1 \leq m \leq r-1$.
		Assume now that $\tau(q) < \frac{1}{m}$ for all $q \in \M$.
		Then the phase map $P:\B\to \B$ is $\C^m$. It additionally follows that $(\B,P,\M)$ is a $\C^{m}$ fibered manifold with $(n-k)$-dimensional Euclidean fibers.
				
		\item $\M$ has unique asymptotic phase $P:\B\to\B$. 		
    \end{enumerate}	
    	
	\end{Prop-hand}

	Proposition \ref{prop:phase} says that if $\nu(p)<1$ and $\sigma(p)<1$, for every $p\in\M$, then $\M$ has unique asymptotic phase.
	It further says that if $\tau(p)<\frac{1}{m}, 1\leq m \leq r-1$ also holds for every $q\in\M$, then the phase map $P$ is $\C^m$.
	
	In stating the next Proposition, we need the following definition.
	
	\begin{Def}\label{def:C1_close}
		Let $g:M\to N\subseteq \R^n$ and $f:M\to N\subseteq\R^n$ be two $\C^1$ maps, where $M$ is a $\C^1$ manifold and $N$ is a $\C^1$ submanifold of $\R^n$. 
		Let $\theta > 0$.
		We say that $g$ and $f$ are \concept{ $\C^1$ $\theta$-close} if:
		\begin{align*}
		\sup_{x\in\B}\|g(x)-f(x)\|&< \theta\\
		\sup_{x\in\B}\|\D g_x-\D f_x\|  &< \theta.
		\end{align*}	
		We will sometimes say that two $\C^1$ maps $f$ and $g$ are \concept{$\C^1$-close} to mean that $f$ and $g$ are $\C^1$ $\theta$-close for some $\theta$ sufficiently small for the present context.
		If $f$ and $g$ are $\C^1$ $\theta$-close, we will sometimes refer to $g$ as a \concept{$\C^1$-small} perturbation of $f$.
		Given two embedded submanifolds $M_1,M_2\subseteq \R^n$, we say that $M_1$ and $M_2$ are \concept{$\C^1$ $\theta$-close} if there exist $\C^1$ embeddings $f_i:M\to \R^n$, $i=1,2$, with each $f_i$ a $\C^1$ diffeomorphism onto $M_i$, such that $f_1$ and $f_2$ are $\C^1$ $\theta$-close.
		Similarly to the case of maps, we will also sometimes simply say $M_1$ and $M_2$ are \concept{$\C^1$-close}.
	\end{Def}
	
	The following restatement of Proposition \ref{prop:persistence} in \S \ref{sec:motivation} is a combination of results from \citet{fenichel1971persistence,fenichel1973asymptotic,fenichel1977asymptotic}, and Theorem 4.1 of \citet{hirsch1977}.

\begin{Prop-hand}{2} %\label{prop:persistence}
	Let $\M$ be a compact $\C^r$ $k$-dimensional embedded submanifold of $\R^n$, invariant under the flow $\phi_t(\cdot)$ of the vector field $f:\Q\to \T \Q$ defined on an open neighborhood $\Q \subset \R^n$ of $\M$.
	Assume that for all $p \in \M$, $\nu(p) < 1$ and $\sigma(p)<\frac{1}{r}$.
	Then for $\theta$ sufficiently small, the following holds:
	\begin{enumerate}
		\item Let $g:\Q \to \T \Q$ be another $\C^r$ vector field which is sufficiently $\C^1$ $\theta$-close to $f$.
		Then there is a unique $\C^r$ embedded submanifold $\M'$, $\C^r$ diffeomorphic to $\M$, $\C^1$-close to $\M$, and invariant under the flow of $g$.
		Furthermore, the fibers $W_q$ persist; i.e., there is a unique invariant fibering of the stability basin of $\M'$ by $\C^r$ manifolds $W_{q'}'$ satisfying all of the properties with respect to $g$ and $\M'$ which were satisfied by the manifolds $W_q$ with respect to $f$ and $\M$.
		The fibers of $W_{q'}'$ are $\C^1$-close to those of $W_q$ on $\B \cap \B'$.
        		$\M'$ has unique asymptotic phase $P':\B'\to\B'$ whose fibers are $W_{q'}'$, and $P':\B' \to \B'$ is a continuous function.
		
		\item Let $1 \leq m \leq r-1$.
		Assume now that $\tau(q) < \frac{1}{m}$ for all $q \in \M$.
		Then the phase map $P':\B'\to \B'$ is also $\C^k$ if $g$ is sufficiently $\C^1$-close to $f$, and also fit together to form a $\C^m$ fibered manifold $(\B',P',\M')$ with $(n-k)$-dimensional Euclidean fibers.
		Under these conditions, $\M'$ has unique $\C^m$ asymptotic phase.
		
	\end{enumerate}
\end{Prop-hand}

Proposition \ref{prop:persistence} is a robustness result; it gives conditions under which $\M$ and its unique $\C^m$ asymptotic phase persist under $\C^1$-small perturbations by $\C^r$ vector fields.

We restate Proposition \ref{prop:converse_mane} from \S \ref{sec:motivation} here.
Proposition \ref{prop:converse_mane} is due to \citet{mane1978persistent}.

\begin{Prop-hand}{3} %\label{prop:converse_mane}
	Let $\M$ be a compact $\C^1$ invariant manifold of the $\C^1$ vector field $f$ which persists under $\C^1$-small perturbations to $f$.
	Then $\M$ is normally hyperbolic.
\end{Prop-hand}

\section{Proofs of \S \ref{sec:our_persistence} results}\label{app:proofs_persistence}

We have shown that under the flow induced by the vector field $f$ on $\M$, $\M$ is asymptotically stable with basin of attraction equal to $\B$.
We have also shown that $\M$ is exponentially stable on a neighborhood $U_E \supset \M$, with exponential rate $\mu$ proportional to $\min_{x \in \bar U_E}\alpha(x)$. 
We now show that if $\min_{x \in \bar U_E}\alpha(x)$ (and hence $\mu$) is chosen sufficiently large, $\M$ can be made $k$-normally hyperbolic for any $k \in \N$.

Note that notation in this section such as $A_t, B_t, \Pi^{\Nor}, \nu, \sigma,$ and $\tau$ is defined in Appendix \ref{app:NHIMs}.

\begin{Lem}\label{lem:Dphi_bound_invariant}
	Let $\mu = \frac{k_4 k_3}{2 k_2}$ be as in Proposition \ref{prop:exp_stable}, where $k_4 = \min_{x \in \bar U_E}\alpha(x)$. 
	There exists $\bar K > 0$ such that for all $t > 0$ and all $p \in \M$:
	\begin{equation*}
	\|\D \phi_t(p)|_{\ker \D P_p}\| < \bar K e^{-\mu t}.
	\end{equation*}
\end{Lem}

\begin{proof}
	The fibers of $P$ are $\C^r$ manifolds transverse to $\M$ depending continuously (actually, in a $\C^r$ manner) on their basepoint in $\M$. 
	Since $\M$ is compact, it follows that (shrinking $U_E$ if necessary) there exists $L > 0$ such that for any $p \in \M$ and $q \in P^{-1}(p)$, 
	\begin{equation}
	\|q-p\| \leq L d_q,
	\end{equation}
	where $d_q$ is the distance from $q$ to $\M$.

	Let $p \in \M$, $v \in \ker \D P_p$.
	Identifying $\T_p \B$ with $\R^n$ and using the fact that $P^{-1}(p)$ is a $\C^r$ manifold with $v \in \T_p P^{-1}(p)$, there exists $q \in P^{-1}(p)$ such that 
	\begin{equation}
	q = p + v + \bo(\|v\|^2). 
	\end{equation}
	This fact together with the theorem on differentiability of flows (\citet{odeHirschSmale} page 299) shows that
	\begin{equation*}
	\|\phi_t(q) - \phi_t(p) - \D(\phi_t)_p v\| \leq \bo(\|v\|^2),
	\end{equation*}
	from which it follows that
	\begin{equation}
	\|\D(\phi_t)_p v\| \leq \|\phi_t(q)-\phi_t(p)\| + \bo(\|v\|^2).
	\end{equation}
	The continuity of ODE solutions  with respect to initial conditions estimate (\citet{odeHirschSmale} page 169) shows that for any fixed $t$, $\bo(\|v\|^2) = \bo(\|\phi_t(q)-\phi_t(p)\|^2)$.
	Since $U_E$ is a region of exponential stability and the fibers of $P$ are invariant under the flow, we have $\|\phi_t(p)-\phi_t(q)\|\to 0$ as $t\to\infty$ and hence this estimate actually holds uniformly for $t>0$.
	It follows that if we take $U_E$ (and hence $\|v\| = \|q-p\|$) sufficiently small, we have for all $t>0$:
	\begin{equation}
	\bo(\|\phi_t(q)-\phi_t(p)\|^2) \leq A \|\phi_t(q)-\phi_t(p)\| \leq AL d_{\phi_t(q)} \leq ALCd_qe^{-\mu t},
	\end{equation}
    for some $A > 0$ and where $C$ is as in the proof of Proposition \ref{prop:exp_stable}.
	This fact, together with the invariance of the fibers of $P$ under the flow and the result of Proposition \ref{prop:exp_stable}, shows that
	\begin{equation}
	\|\D (\phi_t)_p v\| \leq L d_{\phi_t(q)} + \bo(\|v\|^2) \leq [LC+BALC]d_qe^{-\mu t} \leq [LC+BALC]\|v\|e^{-\mu t},
	\end{equation}
	since $\|v\| \leq d_q$, for some $B >0$.
	This completes the proof with $\bar K = LC + BALC$.
\end{proof}

\begin{Co}\label{co:Dphi_bound_NM}
	Let $\mu = \frac{k_4 k_3}{2 k_2 }$ be as in Proposition \ref{prop:exp_stable}, where $k_4 = \min_{x \in \bar U_E}\alpha(x)$.
	Let $\Nor \M$ be the normal bundle of $\M$, and let $\Pi^{\Nor}:\T\R^n|_\M \to \Nor \M$ be the orthogonal projection defined in Appendix \ref{app:NHIMs}.
	Then there exists $K > 0$ such that  for all $t > 0$ and all $p \in \M$:
	\begin{equation*}
	\|\Pi^{\Nor} \D\phi_t(p)|_{\Nor_p\M}\| < K e^{-\mu t}.
	\end{equation*}
	
\end{Co}
\begin{proof}
	$\Ver \B$ and $\Nor \M$ are $\C^{r-1}$ vector bundles over $\M$ of equal dimension and each transverse to $\T \M$.
	It follows that we have a $\C^{r-1}$ linear operator-valued map $\Gamma:\M \to L(\Ver \B,\T \M)$ such that that $\Nor \B$ can be identified with the \concept{graph} of $\Gamma $; precisely, for each $p \in \M:$ $\Gamma(p)$ is a linear map $\Gamma(p):\Ver_p\B \to \T_p\M$ such that  $\Nor_p \B = \{v + \Gamma(p)v|v\in\Ver_p\B\}$.
	Compactness of $\M$ implies $\|\Gamma\|:=\max_{p\in\M}\|\Gamma(p)\|<\infty$ exists.
	
	Let $w \in \Nor_p \M$. 
	The preceding paragraph implies that $w = v + \Gamma(p)v$ for a unique $v \in \ker \D P_p = \Ver_p\B$.
	It follows that $\|w\| \leq (1+\|\Gamma\|)\|v\|$.
	Since $\T\M$ is invariant under $\D\phi_t$, $\Pi^N \D\phi_t \Gamma(p)v = 0$ and hence $\|\Pi^{\Nor} \D\phi_t(p) w\| = \|\Pi^{\Nor} \D\phi_t(p) v\| \leq \|\D\phi_t(p)v\| = \|\D\phi_t(p)|_{\ker \D P_p}v\| \leq \|\D \phi_t(p)|_{\ker \D P_p}\| \|v\| \leq (1+\|\Gamma\|) \|\D \phi_t(p)|_{\ker \D P_p}\|\|w\|$. It follows that $\|\Pi^{\Nor} \D\phi_t(p)|_{\Nor_p\M}\| < K e^{-\mu t}$ with $K:= (1+\|\Gamma\|)\bar K$, where $\bar K$ is as in Lemma \ref{lem:Dphi_bound_invariant}.
\end{proof}

\begin{Lem}
	Define $L:= \max_{p \in \M}\|\D f(p)\| = \max_{p\in\M}\|\D f_0(p)\|$.
	Then
	\begin{align*}
	\|\D(\phi_t)_p|_{\T_p\M}\| \leq \|\D(\phi_t)_p\| < e^{L |t|}.
	\end{align*}
\end{Lem}

\begin{proof}
	Note that the first inequality above is trivial since restricting a linear operator always decreases its norm.
	Now for any $p \in \M$, $\D \phi_t$ satisfies the equation (\citet{odeHirschSmale} pages 300-302):$$\frac{\partial}{\partial t}\D\phi_t(p) = \D f(\phi_t(p))\D\phi_t(p).$$
	Integrating this equation, taking norms, and using the fact that $\D\phi_0(p)$ is the identity shows that 
	\begin{align*}
	\|\D \phi_t(p)\| &\leq 1 + \int_{0}^{t}\|\D f(\phi_s(p))\|\|\D\phi_s(p)\|\,ds\\
	&\leq 1 + L \int_{0}^{t}\|\D\phi_s(p)\|\,ds,
	\end{align*}
	so it follows from Gr\"onwall's Lemma (\citet{odeHirschSmale} page 169) that
	$\|\D\phi_t(p)\| \leq e^{Lt}$.
	Reversing time and repeating the above analysis shows that $$\|\D\phi_t(p)\|\leq e^{L|t|}.$$
	
\end{proof}
\begin{Lem}\label{th:NHIM_rate_conditions}
	Let $\mu = \frac{k_4 k_3}{2 k_2}$ be as in Proposition \ref{prop:exp_stable}, where $k_4 = \min_{x \in \bar U_E}\alpha(x)$,
	and choose $\alpha:\B\to \R$ so that 
	$$k_4 > r \frac{2 k_2}{k_1}L.$$
	Then for every $p \in \M:$
	\begin{align*}
	\nu(p) < 1 \qquad \sigma(p) < \frac{1}{r} \qquad \tau(p) < \frac{1}{r-1}.
	\end{align*}
	
\end{Lem}

\begin{proof}
	It follows immediately from Corollary \ref{co:Dphi_bound_NM} that $\nu(p) < e^{-\mu} < 1$ for every $p \in \M$.
	
	Next, let any $s \in \R$. 
	Our assumption on $k_4$ implies $\mu > r L$.
	We have
	\begin{align*}
	\|A_t(p)\|\|B_t(p)\|^s &= \|\D\phi_{-t}(p)|_{\T_p\M}\|\|\Pi^{\Nor} \D \phi_{t}(\phi_{-t}(p))|_{\Nor_{\phi_{-t}(p)}\M}\|^s \\
	&\leq A K e^{L t}e^{-s \mu t}\\
	&= AK e^{L(1 - sr)t},
	\end{align*}
	which tends to zero for any $s > \frac{1}{r}$ for any $p \in \M$, which shows that $\sigma(p) < \frac{1}{r}$ for any $p \in \M$.
	Similarly,
	\begin{align*}
	\|A_{-t}(p)\| \left(\| B_t(\phi_t(p))\| A_t(\phi_t(p))\|\right)^s & =  \|\D\phi_{t}(p)|_{\T_p\M}\| \left( \|\Pi^{\Nor} \D \phi_{t}(p)|_{\Nor_{\phi_{-t}(p)}\M}\| \|\D\phi_{-t}(\phi_t(p))|_{\T_p\M}\|\right)^s\\
	& \leq A^{1+s}K^se^{[L(1+s)-s\mu]t}\\
	& \leq A^{1+s}K^s e^{[L(1+s)-srL]t}\\
	&= A^{1+s}K^s e^{L(1+s-sr)t},
	\end{align*}
	which tends to zero if $1+s(1-r)< 0$, or for any $s$ such that $s > \frac{1}{r-1}$.
	It follows that for every $p \in \M$, $\tau(p) < \frac{1}{r-1}$. 
\end{proof}

The following theorem is a more precise restatement of Theorem \ref{th:persistence_of_our} in \S \ref{sec:our_persistence}.

\begin{Th-hand}{2}
	Assume $r > 3$.
	Let $\mu = \frac{k_4 k_3}{2 k_2}$ be as in Proposition \ref{prop:exp_stable}, where $k_4 = \min_{x \in \bar U_E}\alpha(x)$,
	and choose $\alpha:\B\to \R$ so that 
	$$k_4 > r \frac{2 k_2}{k_1}L.$$
	Then there exists $\theta > 0$ sufficiently small such that if $g:\B \to \T\B$ is another $\C^{r-1}$ vector field such that
	\begin{align*}
	\sup_{x\in\B}\|g(x)-f(x)\| &< \theta\\
	\sup_{x\in\B}\|\D g(x)-\D f(x)\|  &< \theta,
	\end{align*}	
	then there exists an open set $\B^g\subseteq \B$ and a $\C^{r-1}$ exponentially stable normally hyperbolic submanifold $\M^g$ $\C^{r-1}$ diffeomorphic to $\M$ and $\C^1$ close to $\M$.
	The stability basin of $\M^g$ is $\B^g$.
	$\M^g$ has the unique asymptotic phase property with a $\C^{r-2}$ phase map $P^g:\B^g\to\M^g$ making $(\B^g,P^g,\M^g)$ into a $\C^{r-2}$ fibered manifold with $(n-k)$-dimensional Euclidean fibers.
	The fibers of $P^g$ are $\C^1$ close to the fibers of $P$ on $\B^g$. 	
\end{Th-hand}
\begin{proof}
	This is an immediate corollary of Lemma \ref{th:NHIM_rate_conditions} and Propositions \ref{prop:phase} and \ref{prop:persistence}.
\end{proof}

\bibliographystyle{plainnat}
\bibliography{references}
\end{document}